\documentstyle[12pt] {article}
\def\@height{height}
\def\@depth{depth}
\def\@width{width}
\newtheorem{thm}{Theorem}[section]
\newtheorem{lem}[thm]{Lemma}

\newtheorem{prop}[thm]{Proposition}

\def\an{{< \! {\! {\! {\! {\! )} }}}\, }}
\def\a{{\alpha}}
\def\B{{\beta}}
\def\g{{\gamma}}
\def\G{{\Gamma}}
\def\de{{\delta}}
\def\D{{\Delta}}
\def\si{{\sigma}}
\def\la{{\lambda}}
\def\ka{{\kappa}}
\def\f{{\varphi}}
\def\n{{\nabla}}
\def\o{{\omega}}
\def\O{{\Omega}}

\def\vt{{ {\hat \vartheta} }}
\def\e{{\varepsilon}}

\def\p{{\partial}}
\def\T{{\hspace {0.2cm}}}
\def\TT{{\hspace {0.5cm}}}
\def\nk{{ {\nabla}_{\! \! \hspace{0.3mm} \xi }
\! \hspace{0.8mm} }} 
\def\1{{\rm 1 \!\hspace{0.1mm}\!I}}
\def\bbbr{{\rm I\!R}} 
\def\bbbn{{\rm I\!N}} 

\def\bbbc{{\mathchoice      {\setbox0=\hbox{$\displaystyle\rm
C$}\hbox{\hbox to0pt{\kern0.4\wd0\vrule height0.9\ht0\hss}\box0}}
{\setbox0=\hbox{$\textstyle\rm C$}\hbox{\hbox
to0pt{\kern0.4\wd0\vrule height0.9\ht0\hss}\box0}}
{\setbox0=\hbox{$\scriptstyle\rm C$}\hbox{\hbox
to0pt{\kern0.4\wd0\vrule height0.9\ht0\hss}\box0}}
{\setbox0=\hbox{$\scriptscriptstyle\rm C$}\hbox{\hbox
to0pt{\kern0.4\wd0\vrule height0.9\ht0\hss}\box0}}}}

\def\bbbq{{\mathchoice               {\setbox0=\hbox{$\displaystyle\rm
Q$}\hbox{\raise
0.15\ht0\hbox to0pt{\kern0.4\wd0\vrule height0.8\ht0\hss}\box0}}
{\setbox0=\hbox{$\textstyle\rm Q$}\hbox{\raise
0.15\ht0\hbox to0pt{\kern0.4\wd0\vrule height0.8\ht0\hss}\box0}}
{\setbox0=\hbox{$\scriptstyle\rm Q$}\hbox{\raise
0.15\ht0\hbox to0pt{\kern0.4\wd0\vrule height0.7\ht0\hss}\box0}}
{\setbox0=\hbox{$\scriptscriptstyle\rm Q$}\hbox{\raise
0.15\ht0\hbox to0pt{\kern0.4\wd0\vrule height0.7\ht0\hss}\box0}}}}

\def\bbbz{{\mathchoice {\hbox{$\textstyle\sf Z\kern-0.4em Z$}}
{\hbox{$\textstyle\sf Z\kern-0.4em Z$}}
{\hbox{$\scriptstyle\sf Z\kern-0.3em Z$}}
{\hbox{$\scriptscriptstyle\sf Z\kern-0.2em Z$}}}}

\def\R{\bbbr}
\def\N{\bbbn}
\def\Z{\bbbz}
\def\C{\bbbc}

\title{Sharp Semiclassical Estimates  
for the Number of Eigenvalues 
Below a Degenerate Critical Level }
\date{}
\begin{document}

\setcounter{section}{0} 

\maketitle
 
{Lech Zielinski} \par 
 {\footnotesize {
 {\it LMPA, Centre Mi-Voix, Universit\'e du Littoral,  
B.P. 699, 62228 Calais Cedex, France  \par 
e-mail Lech.Zielinski@lmpa.univ-littoral.fr and IMJ, Math\'ematiques, 
case 7012, \par 
 Universit\'e Paris 7, 2 Place Jussieu, 75251 Paris Cedex 05, France }  }}

\par \smallskip 
{\footnotesize {
{\bf Abstract:} We consider the semiclassical asymptotic behaviour of 
the number of eigenvalues smaller than $E$ for elliptic operators in 
$L^2({\R}^d)$. We describe a method of obtaining remainder estimates   
related to the volume of the region of the phase space in which 
the principal symbol takes values belonging to the intervals $[E';\; E'+h]$, 
where $E'$ is close to $E$. This method allows us to derive 
sharp remainder estimates $O(h^{1-d})$ for a class of symbols 
with critical points and non-smooth coefficients.
\par \smallskip 
{\bf Mathematics Subject Classification (2000):} 35P20
\par \smallskip 
{\bf Keywords:} semiclassical approximation, eigenvalue asymptotics, 
critical energy
 } }

\section{Introduction}

We assume that for $h\! \in ]0;\; h_0]$ the differential 
operators $A_h=a(x,hD,h)$ are self-adjoint in $L^2({\R}^d)$ and  
the symbol $a(x,\xi ,h)=\sum_{0\le n\le N}h^na_n(x,\xi )$ 
is sufficiently regular. If $E\in \R$ satisfies 
 $$
E< \liminf_{|x|+|\xi |\to \infty } a_0(x,\xi ) \eqno (1.1)
$$
and $h_0$ is small enough,  then the spectrum of $A_h$ is discrete 
in $]-\infty ;\; E]$ 
and it is natural to ask whether the counting function 
${\cal N}(A_h,E)$ (i.e. the number of 
eigenvalues smaller than $E$ counted with their multiplicities) 
satisfies the semiclassical asymptotic formula 
 $$
 {\cal N}(A_h,E)=(2\pi h)^{-d}c_E +O(h^{\mu -d})\T {\rm as }\T h\to 0, 
\eqno (1.2)
 $$   
where $\mu >0$ and
 $$
c_E=\int_{a_0(x,\xi )<E}dxd\xi =
{\rm vol}\, \{ (x,\xi )\in {\R}^{2d}:\; a_0(x,\xi )<E\} .\eqno (1.2')
 $$
The most powerful approach of studying semiclassical asymptotics 
has its origin in the microlocal analysis of L. H\"ormander [12]. 
Since the first papers of J. Chazarain [7] and B. Helffer, D. Robert [10],  
this approach has been used in numerous works, cf. 
the monographs [9], [15], [21]. A basic result says that 
(1.2) holds with $\mu =1$ if $E$ is not a critical value of $a_0$ [i.e. 
 $a_0(x,\xi )=E\, \Rightarrow \n a_0(x,\xi )\ne 0$] 
and we refer to the papers T. Paul, A. Uribe [20] 
and M. Combescure, J. Ralston, D. Robert [8], giving  
more precise estimates in relation with the periodic orbits 
of the Hamiltonian flow of $a_0$. 
\par \smallskip 
In this paper we investigate the case when the critical set 
 $$
{\cal C}^{a_0}_E=\{ (x,\xi )\in {\R}^{2d}:\; a_0(x,\xi )=E\hbox{ and }
 \n a_0(x,\xi )=0\} \eqno (1.3)
 $$
is not empty and we consider elliptic operators with 
non-smooth coefficients. Below we enumerate different methods  
and works treating this problem. 
\par 
\smallskip 
-- The analysis of oscillatory integrals.  
If ${\cal C}^{a_0}_E$ is a smooth manifold and the Hessian matrix of $a_0$ 
is transversely non-degenerate, then the semiclassical spectral asymptotics 
can be obtained from the analysis 
described in the paper of R. Brummelhuis, T. Paul, A. Uribe [2]. 
This approach was developed to study the contribution of periodic orbits 
under some geometrical assumptions on the flow (cf.  
D. Khuat-Duy [18], B. Camus [3, 4]) and recent results of B. Camus [5, 6] 
concern the case of a totally degenerate critical point of $a_0$. 
The oscillatory integrals being degenerate in the case of 
a degenerate Hessian matrix, the principal difficulty of this approach 
appears in suitable generalizations of the stationary phase method.   
\par 
\smallskip 
-- The multiscale analysis developed by V. Ivrii [15].  
This method was extended to treat elliptic operators 
with non-smooth coefficients in the paper V. Ivrii [17] 
(cf. also V. Ivrii [16] and M. Bronstein, V. Ivrii [1]).  
\par 
\smallskip 
-- The approximative spectral projector method 
of M. A. Shubin, V. A. Tulovskii [22]. 
The application of this method  
 to our problem was described in the monograph of S. Z. Levendorskii [19]  
and it gives remainder estimates involving the volume of 
suitable regions determined by $a_0$ in the phase space 
valid without any additional assumptions on the Hessian matrix of $a_0$.
After the improvement of L. H\"ormander [13], 
for evry $\mu <{2\over 3}$ one can find  a constant 
$C_{\mu}>0$ such that for $h\in ]0;\; h_0]$ one has  
 $$
\left| { {\cal N}(A_h,E)-(2\pi h)^{-d}c_E } \right| \le C_{\mu}\, h^{-d}
{\cal R}^{a_0}_E(h^{\mu}) \eqno (1.4)
 $$
where  
 $$
 {\cal R}^{a_0}_E(h^{\mu})={\rm vol}\, \{ (x,\xi )\in {\R}^{2d}:\; 
|a_0(x,\xi )-E|\le h^{\mu} \} . 
 $$
\par \smallskip 
-- The method of integrations by parts used in [23, 26]. This method 
allowed us to show that the estimates (1.4) 
are still valid for $\mu <1$. 
\par \smallskip 
\T In this paper we show how to generalize the method of 
[26] to recover estimates with $\mu =1$. Our aim is to show that 
for every $\e >0$ it is possible to find a constant 
$C_{\e}>0$ such that for $h\in ]0;\; h_0]$ one has  
 $$
\left| { {\cal N}(A_h,E)-(2\pi h)^{-d}c_E } \right| \le C_{\e}\, h^{-d}
 {\cal R}^{\e ,a_0}_E(h), \eqno (1.5)
 $$ 
where 
 $$
 {\cal R}^{\e ,a_0}_E(h)=h+\sup_{E'\in [E-h^{1-\e};\; E+h^{1-\e}]} 
 {\rm vol}\, \{ (x,\xi )\in {\R}^{2d}:\; 
|a_0(x,\xi )-E'|\le h \} .
 $$

Using a regularization procedure similarly as in [24, 26] we can 
show that these estimates are valid for elliptic operators with coefficients 
which have second order derivatives H\"older continuous. \par 
It is easy to see that one can always find constants $C,c>0$ such that 
 $$
 {\rm vol}\, \{ (x,\xi )\in {\R}^{2d}:\; 
|a_0(x,\xi )-E'|\le h \}  \le Ch^c \eqno (1.6)
 $$    
and the asymptotic formula (1.2) holds with $\mu =1$ 
if additional properties of $a_0$ ensure the estimate (1.6) with $c=1$ 
for $E'\in [E-h^{1-\e};\; E+h^{1-\e}]$.  
\par 
\T The main part of this paper is devoted to the proof of 
a microlocal trace formula in the region 
 $$
\{ (x,\xi )\in {\R}^{2d}:\; |\n a_0(x,\xi )|\ge h^{{\de}_0} \} \eqno (1.7)
 $$ 
where ${\de}_0\in ]0; 1/2[$. This result allows us to derive 
the asymptotic formula (1.5) under the assumption that 
the Hessian matrix of $a_0$ is of rank $\ge 2$. 
Indeed, it is easy to see (cf. Section 6) that the last assumption 
ensures the fact that the volume of   
 $$
 \{ (x,\xi )\in {\R}^{2d}:\; |\n a_0(x,\xi )|\le h^{{\de}_0} \} \eqno (1.7') 
 $$ 
is $o(h)$ if ${1\over 2}-{\de}_0$ is small enough and the corresponding 
contribution can be included in the right hand side of (1.5) due to (1.6). 
 The next paper [27] will present a trace formula in the region (1.7$'$),
 completing the proof of (1.5) without any hypotheses on the Hessian of $a_0$.
\par \bigskip

 Assume $0<r_0<1$. We write $a\in C_b^{2+r_0}({\R}^d)$   
if and only if the function $a:{\R}^d\to \C$ satisfies the 
following conditions
 $$
 {\p}^{\a}a\in L^{\infty}({\R}^d)\T {\rm if}\T |\a |\le 2,\eqno (1.8) 
 $$
 $$
 |{\p}^{\a}a(x)-{\p}^{\a}a(y)|\le C|x-y|^{r_0}\T {\rm if}\T |\a |=2, \T 
 x,y\in {\R}^d.\eqno (1.8')
 $$
Let $m_0\in {\N}$ and 
for $\nu, {\bar \nu}\in {\N}^d$, $|\nu |, |{\bar \nu}|\le m_0$ let   
$a_{\nu, {\bar \nu}}=a_{{\bar \nu},\nu}\in 
C_{\rm b}^{2+r_0}({\R}^d)$ be real-valued and such that 
 $$
 \sum_{|\nu |=|{\bar \nu}|=m_0}
 a_{\nu ,{\bar \nu}}(x){\xi}^{\nu +{\bar \nu}}\ge c_0|\xi |^{2m_0}
 \TT (x,\xi \in {\R}^d)  \eqno (1.9)
 $$
holds for a certain constant $c_0>0$. Let ${\cal A}_h$ be the quadratic form  
 $$
 {\cal A}_h[\f ,\psi ]=\sum_{|\nu |,|{\bar \nu}|\le m_0}
 (a_{\nu ,{\bar \nu}}(hD)^{\nu}\f ,\, (hD)^{\bar \nu}\psi ),\eqno(1.10)
 $$
where  $\f ,\psi \in C_0^{m_0}({\R}^d)$, \T 
$(hD)^{\nu}=(-ih)^{|\nu |}{{\p}^{\nu}\over \p x^{\nu}}$ 
and  $(\cdot ,\cdot )$ is the scalar product of $L^2({\R}^d)$.
 Due to (1.9), ${\cal A}_h$ is 
bounded from below and its closure defines a self-adjoint operator $A_h$. 
Usually $A_h$ is expressed formally as
 $$
 A_h=\sum_{|\nu |,|{\bar \nu}|\le m_0}(hD)^{\bar \nu} \left( {
   a_{\nu ,{\bar \nu}}(x)(hD)^{\nu} } \right)   . \eqno(1.10')
 $$   
 Moreover we denote 
 $$
 a_0(x,\xi )=\sum_{|\nu |,|{\bar \nu} |\le m_0}a_{\nu ,{\bar \nu}}(x)
 {\xi}^{\nu +{\bar \nu}}. \eqno(1.11)
 $$      
Then we have    

\begin{thm} Let $a_{\nu, {\bar \nu}}\in 
C_b^{2+r_0}({\R}^d)$ be such that $(1.9)$ holds and let 
$A_h$ be self-adjoint operators in $L^2({\R}^d)$ 
defined by $(1.10)$. Let $E\in \R$ be such that $(1.1)$ holds 
(with $a_0$ given by $(1.11)$) and let $h_0>0$ be small enough. 
 \par \smallskip   
a) If $h\in ]0;\; h_0]$ then 
the spectrum of $A_h$ is discrete in $]-\infty ;\; E]$. 
\par \smallskip  
b) If the dimension $d\ge 3$, then for every $\e >0$ one can find a constant 
$C_{\e}>0$ such that $(1.5)$ holds for $h\in ]0;\; h_0]$. 
\end{thm} 
In this paper we show   
\begin{thm} Let $A_h$, $a_0$ and $E$ satisfy the assumptions of Theorem $1.1$. 
Assume moreover that the rank of the Hessian matrix of $a_0$ is greater or equal 
$2$ at every point of the critical set ${\cal C}^{a_0}_E$. 
If the dimension $d\ge 2$, then for every $\e >0$ one can find a constant 
$C_{\e}>0$ such that $(1.5)$ holds for $h\in ]0;\; h_0]$.   
\end{thm} 
The proof of Theorem 1.2 presented in this paper  will be used in [27] 
to derive Theorem 1.1.    
 \par  \smallskip  
{\it Remark.} More general behaviour of coefficients can be considered   
for $x$ such that $a_0(x,\xi )\ge E_0>E$ holds for all $\xi \in {\R}^d$.   
 In particular we have assumed that the coefficients $a_{\nu ,{\bar \nu}}$ 
 are bounded for sake of simplicity, but the same results hold for 
 tempered variations models (cf. B. Helffer, D. Robert [11]). us
 \par  \smallskip     
{\it Plan of the proof.} We begin Section 2 by a description of 
the regularization of non-smooth 
coefficients. It allows us to define the operators $P_h$ with 
smooth coefficients and Theorems 1.1, 1.2, can be deduced from a suitable 
microlocal trace formula for $P_h$. The proof of the trace formula is 
based on the analysis of the evolution group ${\rm exp}(itP_h/h)$ and 
its approximation is described in Section 3. 
\par At the beginning of Section 4 we  
apply the integration  by parts to check the correct trace asymptotics 
of the approximation constructed in Section 3. 
 It remains 
to control the difference between ${\rm exp}(itP_h/h)$ and its 
approximation. Our reasoning is devided in two steps. In Section 4 
we observe that for every ${\bar \e}>0$ one can obtain suitable estimates 
for $|t|\le h^{\bar \e}$ similarly as in [26]. In Section 5 
we use a property of the wave front propagation to show that 
the contribution of the region (1.7) is negligible if $|t|\ge h^{\bar \e}$  
and ${\bar \e}+{\de}_0<{1\over 2}$. In Section 6 we complete the proof    
estimating the volume of the region (1.7$'$).

\section{Regularized problem}

{\bf 2.1 Description of smooth operators} 
\par  \smallskip 
\, We assume ${1\over 2+r_0}<{\de}_0<{1\over 2}$. 
Let ${\g}\in C_0^{\infty}({\R}^d)$ be such that  
$\int {\g} (x)\, dx=1$ and $\int x^{\a}{\g} (x)\, dx=0$  
for $\a \in {\N}^d$ satisfying $1\le |\a |\le 2$.  
We introduce $h$-dependent regularization of coefficients   
 $$
 a_{\nu ,{\bar \nu} ,h}(x) = \int_{{\R}^d} 
a_{\nu ,{\bar \nu} }(y){\g}(h^{-{\de}_0}(x-y))\, h^{-d{\de}_0}\, dy 
\eqno(2.1)
 $$
and define formally self-adjoint differential operators
 $$
 P_h^{\pm}=\sum_{|\nu |,|{\bar \nu}|\le m_0}(hD)^{\nu}
\left( { a_{\nu ,{\bar \nu},h}(x)(hD)^{\bar \nu} } \right) \, \pm 
h(I-h^2\Delta )^{m_0} .  \eqno (2.2)
 $$
We write $P_h^{\pm}$ in the standard form 
 $$
 P_h^{\pm}=\sum_{|\nu |\le 2m_0}p_{\nu ,h}^{\pm}(x)(hD)^{\nu} \eqno (2.3)
 $$
and we use the standard notation $P^{\pm}_h=p_h^{\pm}(x,hD)$ with 
 $$
 p_h^{\pm}(x,\xi )=\sum_{|\nu |\le 2m_0}p_{\nu ,h}^{\pm}(x){\xi}^{\nu}. 
 \eqno (2.4)
 $$
In Section 7 we check the following properties :
 
\begin{lem} Let $a_0$, $A_h$ be as in Theorem $1.1$ and $P^{\pm}_h$, $p^{\pm}_h$ 
as above.  \par  
{\rm (a)} The estimates  
 $$
 |{\p}^{\a}_x{\p}^{\B}_{\xi}p^{\pm}_h(x,\xi )|\le C_{\a ,\B}
 (1+h^{(2+r_0-|\a |){\de}_0})(1+|\xi |)^{2m_0-|\B |}  \eqno(2.5)
 $$
hold for every $\a ,\B \in {\N}^d$ and 
 $$ 
 |{\p}^{\a}_x{\p}^{\B}_{\xi}(a_0-p^{\pm}_h)(x,\xi )|\le C_{\a ,\B}
 (h+h^{(2+r_0-|\a |){\de}_0})(1+|\xi |)^{2m_0-|\B |}  \eqno(2.6)
 $$
hold  if $|\a |\le 2$. \par 
{\rm (b)} Let $h_0>0$ be small enough and consider $h\in  ]0;\; h_0]$.
 Then $A_h^{\pm}$ and the self-adjoint 
realizations of $P_h^{\pm}$ have discrete spectrum  
in $]-\infty ;\; E]$. Moreover the inequalities 
 $$
 P_h^-\le A_h\le P^+_h \eqno (2.7)
 $$ 
hold in the sense of quadratic forms (for $h\in  ]0;\; h_0]$). 
\end{lem} 
\T We deduce (1.5) from suitable 
asymptotics for  $P_h^{\pm}$, observing that (2.7) and 
the min-max principle ensure 
 ${\cal N}(P_h^+,E)\le {\cal N}(A_h,E)\le {\cal N}(P_h^-,E)$. 
Further on we write $P_h$ and $p_h$ instead of $P_h^{\pm}$ and $p_h^{\pm}$. 
 
 \par \bigskip 
{\bf 2.2 Microlocalisation }
\par \smallskip 
 Assume that ${\G}_{\! h}\subset {\R}^{2d}$ for 
 $h\in ]0;\; h_0]$ and denote $\G =({\G}_{\! h})_{h\in ]0;\; h_0]}$.\par  
For $m\in \R$ and ${\de}, {\de}_1\in [0;\; 1[$ satisfying 
${\de}+{\de}_1<1$ we define $S^m_{{\de},{\de}_1}(\G )$ writing   
$b\in S^m_{{\de},{\de}_1}(\G )$  
 if and only if $b=(b_h)_{h\in ]0;\; h_0]}$ 
 is a family of smooth functions $b_h: {\R}^{2d}\to \C$ 
 such that the estimates   
  $$
 \sup_{(x,\xi )\in {\G}_{\! h}} |{\p}^{\a}_{\xi}{\p}^{\B}_xb_h(x,\xi )|
 \le C_{\a ,\B}h^{-m-|\a |{\de}-|\B |{\de}_1}  \eqno (2.8)
 $$
hold for all $\a ,\B \in {\N}^d$. In the case ${\de}={\de}_1$ we abbreviate 
$S^m_{{\de},{\de}}(\G )=S^m_{{\de}}(\G )$. 
\par \smallskip 
  If $b=(b_h)_{h\in ]0;\; h_0]}\in S^m_{{\de},{\de}_1}({\R}^{2d})$  
(i.e. ${\G}_{\! h}={\R}^{2d}$ for all $h$),  then writing 
 $$
 (B_h\f )(x)=(2\pi h)^{-d}\int_{{\R}^d} d\xi \, 
{\rm e}^{ix\xi /h}b_h(x,\xi ) \int_{{\R}^d} dy\; {\rm e}^{-iy\xi /h} \f (y)  
 $$
for $\f \in C_0^{\infty}({\R}^d)$, we define the operators
$B_h=b_h(x,hD)\in B(L^2({\R}^d))$  
satisfying $||b_h(x,hD)||\le Ch^{-m}$ [where $||\cdot ||$ is the norm of 
the algebra of bounded operators $B(L^2({\R}^d))$].
 \par \smallskip 
\T Let ${\1}_Z:\R \to \{ 0,\, 1\}$ be the characteristic 
function of the interval $Z\subset \R$. Then ${\1}_Z(P_h)$
denotes the spectral projector of $P_h$ on $Z$ and 
 $$
 {\cal N}(P_h,E)={\rm tr\,}{\1}_{]-\infty ;\; E]}(P_h). 
 $$  
For a given $t_0>0$ we consider a standard mollifying of ${\1}_Z$ using 
a real valued pair function ${\g}_0\in C_0^{\infty}([-t_0/2;\; t_0/2])$ 
and ${\g}_1={\g}_0\! *\! {\g}_0$ such that ${\g}_1(0)=1$. 
   
The inverse $h$-Fourier transform of ${\g}_1$, given by the formula 
 $$
 {\tilde \g}_h(\zeta )=(2\pi h)^{-1}\int_{\R} dt\;  
{\g}_1(t){\rm e}^{it\zeta /h}\T \T {\rm for}\T \zeta \in \C , \eqno (2.9)
 $$ 
defines a family of holomorphic functions 
satisfying $\int_{\R}{\tilde \g}_h(\la )d\la ={\g}_1(0)=1$ and  
 ${\tilde \g}_h(\la )>0$ for $\la \in \R$. We denote   
   $$
 {\tilde f}^Z_h(\zeta )=\int_Z d\la \, 
{\tilde \g}_h(\zeta -\la ) \eqno (2.9')
 $$
and in Section 7 we show the following result:   
\begin{lem} In order to show the asymptotic formula $(1.5)$ 
it is sufficient to fix $c>0$ and to prove that for every  
$l\in S^0_{0,{\de}_0}({\R}^{2d})$ satisfying 
 $$
{\rm supp\,}l_h\subset \{ v\in {\R}^{2d}: |a_0(v)-E|\le c\} \eqno (2.10)
 $$
 one has the estimates  
 $$
 \left| { {\rm tr\;}{\tilde f}^Z_h(P_h)L_h\,-
 \int_{{\R}^{2d}} {dv \over (2\pi h)^d}\, {\tilde f}^Z_h
 ({\hat p}_h(v ))l_h(v) } \right|
  \le C_{\e}\, h^{-d} \sum_{1\le j\le 2} {\cal R}^{\e ,a_0}_{E_j}(h), \eqno (2.11)
 $$
where $L_h=l_h(x,hD)$, \T  ${\hat p}_h={\rm Re\; }p_h$, \T 
$Z=[E_1;\; E_2]\subset [E-c ;\; E+c]$ and 
the constant $C_{\e}$ is independent of $E_1$, $E_2$.  
\end{lem}  
\par \bigskip 
{\bf 2.3 Partition of the phase space }
\par \smallskip 
Further on ${\bar C},{\bar c}>0$ are constants and we denote 
 $$
 {\hat \G}({\bar C}h^{{\de}_0})=\{ v\in {\R}^{2d}:\, 
 |\nabla a_0(v)|\le {\bar C}h^{{\de}_0} \} , \eqno (2.12)
 $$
 $$
 {\check \G}({\bar c}h^{{\de}_0})=\{ v\in {\R}^{2d}:\, 
 |\nabla a_0(v)|> {\bar c}h^{{\de}_0} \} . \eqno (2.12')
 $$
\par \smallskip 
Using (1.1) we can find   
 ${\overline \G}_{\! 0}$ being a compact subset of ${\R}^{2d}$ such that 
  $$
 \{ v\in {\R}^{2d}: |a_0(v)-E|\le 2c\} \subset {\overline \G}_{\! 0} , \eqno (2.13)  
  $$
where the constant $c>0$ is fixed small enough. 
We will consider $(2.11)$ with $l_h={\hat l}_h+{\check l}_h$, where
 $$
 {\rm supp\;}{\hat l}_h\subset {\hat \G}({\bar C}h^{{\de}_0})\cap 
{\overline \G}_{\! 0},\hspace{5mm} {\rm supp\;}{\check l}_h\subset 
{\check \G}({\bar c}h^{{\de}_0})\cap {\overline \G}_{\! 0} 
 \eqno (2.14)
 $$
and moreover we  introduce an auxiliary cut-off function 
${\tilde l}\in C_0^{\infty}({\R}^d)$ 
such that ${\tilde l}=1$ on ${\overline \G}_{\! 0}$. Further on we denote 
 $$
 {\check L}_h={\check l}_h(x,hD),\hspace{5mm} {\tilde L}_h={\tilde l}_h(x,hD), 
 \eqno (2.15) 
 $$ 
where ${\check l}$, ${\tilde l}$ are as above. Then we have
 $||{\check L}_h(I-{\tilde L}_h^*)||_{\rm tr}=O(h^{\infty})$, 
 where  $||\cdot ||_{\rm tr}$ denotes the trace class norm,  
 $O(h^{\infty})$ means that $O(h^m)$ holds for every $m\in \R$
 and ${\tilde L}_h^*$ denotes the adjoint of ${\tilde L}_h$ in $L^2({\R}^d)$.
Using moreover  the trace cyclicity we obtain  
   $$
 {\rm tr\;}{\check L}_h{\tilde f}^Z_h(P_h)\,=\, 
 {\rm tr\;}{\check L}_h{\tilde f}^Z_h(P_h){\tilde L}_h^*\,+O(h^{\infty}) 
 $$   
and further on we keep the auxiliary cut-off ${\tilde L}_h^*$ to be sure that 
our analysis always concern operators of the trace class. We introduce    
 $$
 f^Z_h(t)=\int_{\R}d\la \; {\rm e}^{-it\la /h}{\tilde f}^Z_h(\la )={\g}_1(t)
 \int_Z d\la \; {\rm e}^{-it\la /h}. \eqno(2.16)
 $$
Then \T  ${\tilde f}_h(\la )=\int_{\R} 
 {dt\over 2\pi h}\; f^Z_h(t)\, {\rm e}^{it\la /h}$    
 \T and consequently
 $$
 {\rm tr\,}{\check L}_h{\tilde f}_h(P_h){\tilde L}_h^*\; =\int_{\R}{dt\over 2\pi h}
 \; f_h(t)\, {\rm tr\,}{\check L}_h{\rm e}^{itP_h/h}{\tilde L}_h^*.  
 \eqno (2.17)
 $$ 
We will prove   
 \begin{prop} Assume that $t_0,\ka >0$ are small enough, $\e >0$ and  
${\check L}_h$, ${\tilde L}_h$ are given by $(2.15)$.  Then 
for every ${\bar N}\in \N$ one can find 
the operators $Q^h_{\bar N}(t)\in B(L^2({\R}^d))$ such that 
$$ 
 \sup_{-t_0\le t\le t_0} |\, {\rm tr\;} \lgroup 
 Q^h_{\bar N}(t)-{\check L}_h{\rm e}^{itP_h/h}){\tilde L}_h^*
\, |\T \le \, C_{\bar N}h^{{\bar N}\kappa -5d-1}  \eqno (2.18)
 $$
and one has the estimates 
 $$
 \left| {\int_{\R}{dt\over 2\pi h}\; f^Z_h(t)\, {\rm tr\,}Q^h_{\bar N}(t)
 {\tilde L}_h^*
 \, -\int_{{\R}^{2d}}{dv \over (2\pi h)^d}\, {\check l}_h(v) 
 {\tilde f}^Z_h({\hat p}_h(v))\, } \right| 
 $$
 $$
  \le C_{\e}\, h^{-d} \sum_{1\le j\le 2} {\cal R}^{\e ,a_0}_{E_j}(h), \eqno (2.19)
 $$
where ${\hat p}_h={\rm Re\; }p_h$, \T 
$Z=[E_1;\; E_2]\subset [E-c ;\; E+c]$ and 
$C_{\e}$ is independent of $E_1$, $E_2$.
 \end{prop} 
Using Proposition 2.3 we find that (2.11) holds if $l_h$, $L_h$ 
are replaced by ${\check l}_h$, ${\check L}_h$. 
Indeed, we have $|f^Z_h(t)|\le C|{\g}_1(t)|$ and 
${\rm supp\;}{\g}_1\subset [-t_0;\; t_0]$, hence we can replace 
${\check L}_h{\rm e}^{itP_h/h}$ by $Q_{{\bar N},h}(t)$ in (2.17) with an error 
$O(h^{{\bar N}\kappa -5d-1})$. 
 \par \smallskip  
\T The construction of $Q^h_{\bar N}(t)$ is presented in Section 3 
and Proposition 2.3 is proved in Sections 4-5.    
At the end of this section 
we introduce classes of symbols describing properties of  
${\check l}_h$ and ${\hat l}_h$. 
 \par \bigskip 

{\bf 2.4 Classes of symbols $S^m_{(0)}$ } 
\par \smallskip 
Further on ${\de}_0$, 
${\overline {\G}}_{\! 0}$ are fixed as before, we fix ${\bar C}>1$, we denote    
 $$
 {\hat \G}_{\! h}={\hat \G}({\bar C}h^{{\de}_0}), \T \T \T \T  \T \T    
 {\check \G}_{\! h}={\R}^{2d}\setminus {\hat \G}_{\! h}=
 {\check \G}({\bar C}h^{{\de}_0}) \eqno (2.20)
 $$ 
and we write 
${\hat \G} =({\hat \G}_{\! h})_{h\in ]0;\; h_0]}$, 
${\check \G} =({\check \G}_{\! h})_{h\in ]0;\; h_0]}$.\par  
\T For $m\in \R$ we define $S^m_{(0)}$ writing 
  $b\in S^m_{(0)}$ if and only if  
 $b=(b_h)_{h\in ]0;\; h_0]}\in 
S^m_{0,{\de}_0}({\check \G})\cap S^m_{{\de}_0}({\hat \G})$ 
and ${\rm supp\,}b\subset {\overline {\G}}_{\! 0}$ (i.e. 
${\rm supp\,}b_h\subset {\overline {\G}}_{\! 0}$ for every $h\in ]0;\; h_0]$). 
\begin{lem} Let $l\in S^0_{0,{\de}_0}$ satisfy  
${\rm supp\,}l\subset {\overline {\G}}_{\! 0}$. Then there exist 
${\check l}, {\hat l}\in S^0_{(0)}$ such that   
$l={\hat l}+{\check l}$, ${\rm supp\,}{\hat l}\subset {\rm supp\,}l$ 
and $(2.14)$ holds if ${\bar c}>0$ is small enough. 
\end{lem} 
{\bf Proof.} Let ${\chi}_0\in C_0^{\infty}(]-1;\; 1[)$ satisfy 
 $0\le {\chi}_0\le 1$ and ${\chi}_0=1$ on $[-{1\over 2};\; {1\over 2}]$. 
We define ${\hat \chi}\in S^0_{{\de}_0}({\overline \G}_{\! 0})$ setting 
 ${\hat \chi}_h(v)=\prod_{|\B |=1}
{\chi}_0({\p}^{\B}{\hat p}_h(v)\cdot 2d/h^{{\de}_0})$. Since 
   $$
 \nabla p_h(x,\xi )=\nabla a_0(x,\xi )+o(h^{{\de}_0})(1+|\xi |^{2m_0}) 
 \eqno (2.21)
 $$
due to (2.6) and $a_0$ is real valued, (2.21) still holds if we replace $p_h$ by 
${\hat p}_h$. However ${\hat \chi}_h(v)\ne 0$ implies    
$|{\p}^{\B}{\hat p}_h(v)|\cdot 2d/h^{{\de}_0}\le 1$ if $|\B |=1$, hence 
we have $|\nabla a_0(v)|\le h^{{\de}_0}(1+o(1))$ 
if moreover $v\in {\overline \G}_{\! 0}$. 
Therefore setting ${\hat l}=l{\hat \chi}$ we obtain 
${\rm supp\;}{\hat l}_h\subset {\hat \G}({\bar C}h^{{\de}_0})$ 
if ${\bar C}>1$ and consequently ${\hat l}\in S^0_{(0)}$. \par 
Next we define ${\check l}=l-{\hat l}\in S^0_{(0)}$. 
If $v\in {\rm supp\;}{\check l}$, then ${\hat \chi}_h(v)\ne 1$ and 
$|{\p}^{\B}p_h(v)|\cdot 2d/h^{{\de}_0}\ge 1/2$ holds for a certain 
$\B \in {\R}^{2d}$ with $|\B |=1$, hence 
$|\nabla p_h(v)|\ge h^{{\de}_0}/(4d^2)$. Therefore      
${\rm supp\;}{\check \chi}_h\subset {\check \G}({\bar c}h^{{\de}_0})$ 
holds if ${\bar c}<1/(4d^2)$. $\triangle$. 


\section{Approximation of the evolution}

{\bf 3.1 Preliminaries} 
\par \smallskip 
We write ${\hat p}_h={\rm Re\;}p_h$ and we consider an approximation of 
${\check L}_h{\rm e}^{itP_h/h}$ in the form   
 $$
 Q^h_{\bar N}(t)=\Bigl( { {\rm e}^{it{\hat p}_h/h} \! \sum_{0\le n\le {\bar N}}
 t^nq^{\circ}_{{\bar N},n,h} }\Bigr) (x,hD),\eqno (3.1) 
 $$
where $q^{\circ}_{{\bar N},n,h}$ will be described in Proposition 3.3.
 We introduce formally  
 $$
  {\tilde Q}^h_{\bar N}(t)=
\hbox{${d\over dt}$}Q^h_{\bar N}(t)-iQ^h_{\bar N}(t)P_h/h  \eqno (3.2)
 $$ 
and require  $Q^h_{\bar N}(0)={\check L}_h$, which allows us to express     
 $$
  Q^h_{\bar N}(t)-{\check L}_h{\rm e}^{itP_h/h}=
 \int_0^t d\tau \, {\tilde Q}^h_{\bar N}(t-\tau ){\rm e}^{i\tau P_h/h}. 
 \eqno (3.3)
 $$
To investigate (3.2) in terms of symbols we introduce the notation  
  $$
  ({\tilde {\cal P}}_{\! \bar N}b)_h(t) ={\rm e}^{-it{\hat p}_h/h} \! \! 
  \left( { \! {\p}_t(b_h(t) {\rm e}^{it{\hat p}_h/h}) -\! \sum_{|\a |\le {\bar N}} 
 \! {h^{|\a |-1}\over \a !\, i^{|\a |-1}}
   {\p}^{\a}_{\xi}  ( b_h(t)\, {\rm e}^{it{\hat p}_h/h}
   {\overline { {\p}^{\a}_xp} }_h ) \! }  \right)  \eqno(3.4)
  $$ 
if $(b_h(t))_{h\in ]0;\; h_0]}\in S^m_{(0)}$ for $t\in \R$.  
 \par \bigskip  
{\bf 3.2 Classes of symbols $S^m_{(N)}$ for $N\in \N$ }
\par \smallskip 
We use the induction with respect to $N\in \N$ in the following definition: \par 
i) if $N=0$ then $S^m_{(N)}=S^m_{(0)}$ was defined at the end of Section 2 ;\par 
ii) we write $(b_h)_{h\in ]0;\; h_0]}\in S^m_{(N+1)}$ if and only if 
 $$
 b_h(x,\xi )=b_{0,h}(x,\xi )+\sum_{|\B |=1}h^{-1/2}b_{\B ,h}(x,\xi )
{\p}^{\B}{\hat p}_h(x,\xi ) \eqno (3.5) 
 $$    
holds with some $(b_{\B ,h})_{h\in ]0;\; h_0]}\in S^m_{(N)}$ for  
  $\B \in {\N}^{2d}$ satisfying  
 $|\B |\le 1$. \par 
It is clear that $S^m_{(N)} \subset S^{m+N/2}_{(0)}$. Moreover  
 $b\in S^m_{(N)}$ if and only if one can write
 $$
 b(x,\xi )=\sum_{ \{ \B \in {\N}^{2d}:\; |\B |\le N\} }
 b_{\B}(x,\xi )h^{-|\B |/2}(\nabla {\hat p}(x,\xi ))^{\B} \eqno (3.6)
 $$
with some symbols $b_{\B}\in S^m_{(0)}$, where 
$\nabla {\hat p}(x,\xi )\in {\R}^{2d}$ and 
 $$
 (t_1,...,t_n)^{({\a}_1,...,{\a}_n)}=\prod_{1\le j\le n}{\! t_j^{{\a}_j}}\T \,   
{\rm for}\T(t_1,...,t_n)\in {\R}^n,\T ({\a}_1,...,{\a}_n)\in {\N}^n. 
 $$  
Using this characterization we find
 $$
 b\in S^m_{(N)},\T {\tilde b}\in S^{\tilde m}_{0,{\de}_0}({\overline \G}_{\! 0})
 \Rightarrow b{\tilde b}\in S^{m+{\tilde m}}_{(N)}, \eqno (3.7)
 $$ 
 $$
 b\in S^m_{(N)},\T {\tilde b}\in S^{\tilde m}_{({\tilde N})}
 \Rightarrow b{\tilde b}\in S^{m+{\tilde m}}_{(N+{\tilde N})}. \eqno (3.8)
 $$  
\begin{lem} {\bf {\rm (a)}} If $b\in S^m_{(N)}$ then   
 ${\p}^{\a}b \in S^{m+|\a |/2}_{(N)}$ for every $\a \in {\N}^{2d}$. \par 
{\bf {\rm (b)}} If $b\in S^m_{(N+1)}$ then 
 ${\p}_{{\xi}_k}b \in S^{m+1/2}_{(N)}$ for $k\in \{ 1,...,d\}$.
\end{lem}  
{\bf Proof (a)} The general statement follows by induction 
with respect to $|\a |$ and we consider only the case $|\a |=1$.   
It is clear that the assertion holds for $N=0$ and reasoning 
by induction we assume that the assertion holds for a given $N\in \N$. 
Let $b\in S^m_{(N+1)}$. Then (3.5) holds with $b_{\B}\in S^m_{(N)}$, hence    
${\p}^{\a}b_{\B}\in S^{m+1/2}_{(N)}$ if $|\a |=1$. Due to (2.5) we have 
${\p}^{\a +\B}{\hat p}\in S^0_{0,{\de}_0}({\overline \G}_{\! 0})$
and writing  
  $$
 {\p}^{\a}(b_{\B}\; h^{-1/2}{\p}^{\B}{\hat p})=
b_{\B}\; h^{-1/2}{\p}^{\a +\B}{\hat p}
+{\p}^{\a}b_{\B}\; h^{-1/2}{\p}^{\B}{\hat p}
 \in S^{m+1/2}_{(N+1)}, \eqno (3.9)
 $$  
we obtain the assertion of Lemma 3.1(a) for $N+1$. 
   \par \smallskip
{\bf (b)} We will show that for $k\in \{ 1,...,d\}$,  
 $|\a |=1$ and $N\in \N$ we have
 $$
 b\in S^m_{(N)}\, \Rightarrow \; {\p}_{{\xi}_k}b\; {\p}^{\a}{\hat p}
\in S^m_{(N)}.  \eqno (3.10(N))
 $$
To begin we consider $N=0$. Since 
${\p}^{\a}{\hat p}\in S^0_{0,{\de}_0}({\overline \G}_{\! 0})$, we have 
 $$
 b\in S^m_{(0)}\Rightarrow {\p}_{{\xi}_k}b\in 
S^m_{0,{\de}_0}({\check \G}) \Rightarrow \; {\p}_{{\xi}_k}b\; {\p}^{\a}{\hat p}
\in S^m_{0,{\de}_0}({\check \G}). \eqno (3.11)
 $$
It is easy to check that  
${\p}^{\a}{\hat p}\in S^{-{\de}_0}_{{\de}_0}({\hat \G}\cap 
{\overline \G}_{\! 0})$ and consequently 
 $$
 b\in S^m_{(0)}\Rightarrow {\p}_{{\xi}_k}b\in S^{m+{\de}_0}_{{\de}_0}
({\hat \G}) \Rightarrow \; {\p}_{{\xi}_k}b\; {\p}^{\a}{\hat p}\in 
S^m_{{\de}_0}({\hat \G}). \eqno (3.12)
 $$
However (3.11) and (3.12) imply 
${\p}_{{\xi}_k}b\; {\p}^{\a}{\hat p}\in S^m_{(0)}$, i.e. (3.10(0)) holds. \par 
\T Reasoning by induction we assume that (3.10($N$)) holds for a given 
$N\in \N$. If $b\in S^m_{(N+1)}$ is given by (3.5), then 
 $$
b_{\B}\in S^m_{(N)}\Rightarrow {\p}_{{\xi}_k}{\p}^{\B}{\hat p}\; b_{\B}
\in S^m_{(N)}\Rightarrow {\p}_{{\xi}_k}{\p}^{\B}{\hat p}\; b_{\B}\; 
h^{-1/2}{\p}^{\a}{\hat p} \in S^m_{(N+1)}. \eqno (3.13)
 $$
Moreover the induction hypothesis ensures
 $$
 b_{\B}\in S^m_{(N)} \Rightarrow {\p}_{{\xi}_k}b_{\B}\; {\p}^{\a}{\hat p}
\in S^m_{(N)} \Rightarrow {\p}_{{\xi}_k}b_{\B}\; {\p}^{\a}{\hat p}\; 
h^{-1/2}{\p}^{\B}{\hat p}\in S^m_{(N+1)}. \eqno (3.14)
 $$
Summing up (3.13) and (3.14) we obtain 
${\p}_{{\xi}_k}(b_{\B}\; h^{-1/2}{\p}^{\B}{\hat p}){\p}^{\a}{\hat p}
\in S^m_{(N+1)}$, which completes the proof of 
(3.10$(N+1)$). 
 \par 
\T If $b\in S^m_{(N+1)}$ is given by (3.5) with $b_{\B}\in S^m_{(N)}$, 
then ${\p}_{{\xi}_k}b_{\B}\, {\p}^{\B}{\hat p}\in S^m_{(N)}$ 
is ensured by (3.10$(N)$) and we obtain the assertion of Lemma 3.1(b) 
writing (3.9) with ${\p}_{{\xi}_k}$ instead of ${\p}^{\a}$. 
 $\triangle$   
\par \bigskip 

{\bf 3.3 Construction of the approximation} 
\par \smallskip  
We define auxiliary classes of symbols ${\check S}^m_{(N)}\subset S^m_{(N)}$
for $N\in \N \setminus \{ 0\}$ as follows:  we write $b\in {\check S}^m_{(N)}$ 
if and only if it is possible to find $b_j\in S^m_{(N-1)}$, 
$j\in \{ 0, 1,...,d\}$, such that   
 $b=b_0+\sum_{1\le j\le d}h^{-1/2}b_j{\p}_{{\xi}_j}{\hat p}$.   
\par 

\begin{lem} Let $b$ be independent of $t$. 
If $b\in S^m_{(N)}$ then 
  $$
 ({\tilde {\cal P}}_{\bar N}b)(t) =\sum_{0\le n\le {\bar N}}t^nb_n \eqno (3.15)
  $$
holds with $b_0\in S^m_{(N)}\subset {\check S}^m_{(N+1)}$ and 
$b_n\in {\check S}^m_{(N+n+1)}$ for $n\in \{ 1,...,{\bar N} \}$.
 \end{lem} 
  {\bf Proof.} To begin we show that $b_0 \in S_{(N)}^m$. We observe that
  $$
  b_0=ih^{-1}({\hat p}-{\overline p})b+ \sum_{1\le |\a |\le {\bar N}} 
{h^{|\a |-1}\over \a !\, i^{|\a |-1}}{\p}^{\a}_{\xi} (b\, 
{\overline {{\p}^{\a}_xp } }\; )  \eqno(3.16)
  $$
and since $p_h(x,hD)$ is self-adjoint, using (2.5) we obtain 
 $$
 |\a |\le 1\Rightarrow {\p}^{\a}({\hat p}-{\overline p})=2i\; {\rm Im\;}{\p}^{\a}p
 \in S^{-1}_{0,{\de}_0}({\overline \G}_{\! 0}). \eqno (3.17)
 $$
We observe that the first term of the 
expresion (3.16) belongs to $S_{(N)}^m$ due to 
(3.17) with $\a =0$. Then $b\in S^m_{(N)}\Rightarrow   
b{\p}_{x_j}{\hat p}\in S^{m-1/2}_{(N+1)}$ and using (3.17) with $|\a |=1$ 
we obtain $b{\overline {{\p}_{x_j}p } }\in S^{m-1/2}_{(N+1)}$. Therefore 
 Lemma 3.1(b) ensures 
${\p}_{{\xi}_j}(b{\overline {{\p}_{x_j}p } })\in S^m_{(N)}$,  
i.e. all terms of (3.16) 
with $|\a |=1$ belong  $S_{(N)}^m$. \par   
 In the next step we consider the terms of (3.16) with $|\a |\ge 2$. \par
Since (2.5) ensures  
$|\a |\ge 2\Rightarrow 
{\p}^{\a}_{x}p\in S^{(|\a |-2){\de}_0}_{0,{\de}_0}({\overline \G}_{\! 0})$, 
we obtain
 $$
 b\; {\overline {{\p}^{\a}_xp } }\in S^{m+(|\a |-2){\de}_0}_{(N)}\Rightarrow 
 h^{|\a |-1}{\p}^{\a}_{\xi}(b\; {\overline {{\p}^{\a}_xp } })\in 
S^{-|\a |+1+|\a |/2+m+(|\a |-2){\de}_0}_{(N)}  
 $$
due to Lemma 3.1(a) and $m+(|\a |-2)({\de}_0-{1\over 2})\le m$ 
gives $b_0 \in S_{(N)}^m$.    
 \par  
 In order to show $b_n\in {\check S}^m_{(N+n+1)}$ for $n\in \{ 2,...,{\bar N} \}$  
we write  
 $$
 b_n=\sum_{ {}^{|\B |\le n\le |\a |\le {\bar N}}_{\hspace{3mm} 
 \B +{\bar \B}\le \a } }
 h^{|\a |-1-n}\; b_{\B ,{\bar \B}}\, (\nk {\hat p})^{\B}\, 
 {\p}_{\xi}^{\bar \B}(b\, {\overline {{\p}^{\a}_x p } }\; ),
 $$
where  $b_{\B ,{\bar \B}}\in S^0_{0,{\de}_0}({\overline \G}_{\! 0})$ 
for $\B ,{\bar \B}\in {\N}^d$. More precisely: in the case $|\B |<n$,  
$b_{\B ,{\bar \B}}$ is a linear combination of terms
 ${\Pi}_{1\le k\le n-|\B |} {\p}_{\xi}^{{\bar \a}(k)}{\hat p}$
 where ${\bar \a}(k)\in {\N}^d$ are such that  
$|{\bar \a}(k)|\ge 2$ for $k\in \{ 1,...,n-|\B |\}$ and 
 $\B +{\bar \B}+\sum_{1\le k\le n-|\B |}{\bar \a}(k)=\a$, implying 
 $$
 |\a |\ge |\B |+|{\bar \B}|+2(n-|\B |)= 
 2n+|{\bar \B}|-|\B |. \eqno (3.18)
 $$
In the case $|\B |=n$ the symbols $b_{\B ,{\bar \B}}$ are constant and (3.18) 
still holds. \par 
Consider first the case $|\a |\ge 2$.  
Then using Lemma 3.1(a) we find 
 $$
 h^{|\a |-1-n}{\p}_{\xi}^{\bar \B}(b\, 
{\overline {{\p}^{\a}_x p } }\; )(\nk {\hat p})^{\B}
\in S_{(N+|\B |)}^{-|\a |+1+n+|{\bar \B}|/2+m+(|\a |-2){\de}_0-|\B |/2} 
 \eqno (3.19)
 $$  
and (3.18) ensures 
 $$
1-2{\de}_0+(2n+|{\bar \B}|-|\B |)/2+|\a |({\de}_0-1)\le 
 (|\a |-2)({\de}_0-1/2)\le 0, \eqno (3.20) 
$$
i.e. all terms corresponding to $|\a |\ge 2$ belong to 
$S^m_{(N+n)}\subset {\check S}^m_{(N+n+1)}$.
\par To complete the proof we observe that in the case $n=|\a |=1$ we have 
 $$
 b\; {\overline {{\p}_{x_j}p } }
\in S^{m-1/2}_{(N+1)}\Rightarrow  
h^{-1}b\; {\overline {{\p}_{x_j}p } }\; {\p}_{{\xi}_j}{\hat p}
\in {\check S}^m_{(N+2)}. \eqno \triangle
 $$
\par \smallskip  
 
  \begin{prop} Let ${\check l}\in S^0_{(0)}$ and ${\bar N}\in \N$. 
Assume that $N\in \{ 0,1,...,{\bar N}\}$. 
Then we can find 
 $$
  q_{{\bar N},N}(t)=\sum_{0\le n\le N}t^nq^{\circ}_{{\bar N},n} 
\eqno (3.21(N)) 
 $$ 
 such that\T $q^{\circ}_{{\bar N},0}={\check l}$, \T 
$q^{\circ}_{{\bar N},1}\in S_{(0)}^0$,  
 $$
  q^{\circ}_{{\bar N},n}\in {\check S}^0_{(n)} \TT
 {\it for}\T n\in \{ 2,...,N\}  \eqno (3.22(N))
  $$
and 
  $$
  {\tilde {\cal P}}_{\bar N}q_{{\bar N},N}(t)=\sum_{N\le n\le N+{\bar N}} 
t^n{\tilde q}^{\, \circ}_{{\bar N},N,n} \eqno (3.23(N))
  $$
holds with  
  $$
  {\tilde q}^{\, \circ}_{{\bar N},N,n}\in {\check S}_{(n+1)}^0 \TT 
{\it for}\T n\in \{ N,...,N+{\bar N} \} .  \eqno (3.24(N))
  $$
\end{prop}
 
{\bf Proof.} If $N=0$ then we take 
$q^{\circ}_{{\bar N},0}={\check l}\in S_{(0)}^0$ 
and Lemma 3.2 with $b={\check l}$ ensures (3.23(0)) and (3.24(0)).  
 Next we assume that the statement of Proposition 3.3 holds for a given
$N\le {\bar N}-1$ and 
using the induction hypothesis $(3.23(N))$ to express
${\tilde {\cal P}}_{\bar N}q_{{\bar N},N}(t)$ we find
  $$
  {\tilde {\cal P}}_{\bar N}q_{{\bar N},N+1}(t)=
  {\tilde {\cal P}}_{\bar N}(t^{N+1}q^{\circ}_{{\bar N},N+1})+
  {\tilde {\cal P}}_{\bar N}q_{{\bar N},N}(t)=
  $$
  $$
  t^N\Bigl( (N+1)q^{\circ}_{{\bar N},N+1}+
{\tilde q}^{\, \circ}_{{\bar N},N,N}\Bigl) 
 \, +\, t^{N+1}{\tilde {\cal P}_{\bar N}}q^{\circ}_{{\bar N},N+1}(t)+
 \sum_{N+1\le n\le N+{\bar N}} t^n\, {\tilde q}^{\, \circ}_{{\bar N},N,n}.
  $$
To obtain $(3.23(N+1))$ we cancel the term with $t^N$ taking 
  $$
  q^{\circ}_{{\bar N},N+1}=-{\tilde q}^{\, \circ}_{{\bar N},N,N}/(N+1), 
  $$
which is an element of ${\check S}^0_{(N+1)}$ by the induction hypothesis 
$(3.24(N))$ 
and (3.24($N+1$)) follows if we develop  
$t^{N+1}{\tilde {\cal P}_{\bar N}}q^{\circ}_{{\bar N},N+1}(t)$ 
as in Lemma 3.2 with $b=q^{\circ}_{{\bar N},N+1}\in S^0_{(N+1)}$.  
Moreover for $N=0$ we have 
$q^{\circ}_{{\bar N},1}=-{\tilde q}^{\, \circ}_{{\bar N},0,0}\in S^0_{(0)}$ 
(due to Lemma 3.2 with $b={\check l}$).  $\triangle$ 
 
\section{Quality of the approximation}

\T This section is devoted to the proof of Proposition 2.3. To begin  
we introduce more notations. We write $q\in {\tilde S}^m_{(0)}$ if and only if  
$q=(q_h)_{h\in ]0;\; h_0]}$ with $q_h\in C^{\infty}({\R}^{3d})$ 
satisfying the estimates   
  $$
 |{\p}^{\a}q_h(x,\xi ,y)|\le C_{\a}
 h^{-m-|\a |{\de}_0} \eqno (4.1)
 $$
for every $\a \in {\N}^{3d}$ and 
 ${\rm supp\, }q_h\subset \, {\overline {\G}}_{\! 0}\times {\R}^d.$  
As before  ${\hat p}_h={\rm Re\;}p_h$ and writing 
 $$
 ({\rm Op}^h_t[q]\f )(x)=\int_{{\R}^{2d}} {dyd\xi \over (2\pi h)^d}\, 
 {\rm e}^{i(x-y)\xi /h+it{\hat p}_h(x,\xi )/h}q_h(x,\xi ,y)\f (y)   
 \eqno (4.2)
 $$ 
for $\f \in C_0^{\infty}({\R}^d)$ we define operators on $L^2({\R}^d)$ 
such that 
 $$
 q\in {\tilde S}^m_{(0)}\, \Rightarrow \, 
 \sup_{-t_0\le t\le t_0}||{\rm Op}^h_t[q]||_{\rm tr}\le Ch^{-m-5d}. \eqno (4.3)
 $$ 
Indeed, (4.3) follows from standard 
estimates of pseudo-differential operators (e.g. [14, Sec. 18]), 
the details are given in the proof of (4.4) in [24].\par  
We observe that   
 $$
 {\rm tr\;}{\rm Op}^h_t[q]\; =\int_{{\R}^{2d}} {dxd\xi \over (2\pi h)^d}\, 
 {\rm e}^{it{\hat p }_h(x,\xi )/h}q_h(x,\xi ,x) \eqno (4.4) 
 $$
and for $b_h\in C_0^{\infty}({\R}^{2d})$ we introduce the notation 
 $$
 J_t^h(b)=\int_{{\R}^{2d}}{dxd\xi \over (2\pi h)^d}\, 
 {\rm e}^{it{\hat p }_h(x,\xi )/h}b_h(x,\xi ). \eqno (4.5)
 $$ 
Using this notation and \T  
$Q^h_{\bar N}(t){\tilde L}_h^*=\sum_{0\le k\le {\bar N}}t^k {\rm Op}_t^h 
[q^{\; \circ}_{{\bar N},k}(x,\xi ){\overline {{\tilde l}(y,\xi )} } ]$ \T 
with ${\tilde l}=1$ on ${\rm supp\;}q^{\; \circ}_{{\bar N},k}$, we find 
the expression 
  $$
 {\rm tr\,}Q^h_{\bar N}(t){\tilde L}_h^*\; = 
\sum_{0\le k\le {\bar N}}t^kJ^h_t(q^{\; \circ}_{{\bar N},k}).  \eqno (4.6)
 $$

\begin{lem} Let $n\in \N$. {\rm (a)} If $b\in S_{(N)}^m$ then 
 $$
 t^nJ^h_t(b)=\sum_{0\le k\le n}t^kJ^h_t(b_{k,n}) \eqno (4.7(n))
 $$ 
holds with some $b_{k,n}\in S_{(\max \{ 0,N-n\} )}^m$ 
for $k\in \{ 0,\dots ,n\}$. \par 
{\rm (b)} If $b\in {\check S}_{(N+1)}^m$, $N\ge 1$,  then $(4.7(n))$ holds with 
$b_{k,n}\in S_{(\max \{ 0,N-n\} )}^m$.
\end{lem}
{\bf Proof.} (a) Reasoning by induction we assume that the statement holds 
for a given $N\in {\N}$. In order to show that the statement still holds for $N+1$ 
instead of $N$ we consider $b\in S_{(N+1)}^m$. Then (3.5) holds with 
$b_{\B}\in S_{(N)}^m$ and  $t^nJ^h_t(b_0)$ 
can be expressed in a suitable way due to the induction hypothesis.  
Then the integration by parts gives 
 $$
 tJ_t^h(h^{-1/2}{\p}^{\B}{\hat p}\; b_{\B})=J_t^h({\tilde b}_{\B})
\T {\rm with}\T {\tilde b}_{\B}=h^{1/2}i{\p}^{\B}b_{\B} \eqno (4.8)
 $$ 
and Lemma 3.1(a) ensures ${\tilde b}_{\B}\in S_{(N)}^m$, i.e.  
(4.8) implies (4.7(1)). Reasoning by induction 
with respect to $n\in \N$ we obtain (4.7($n$)). 
\par (b) If $b\in {\check S}_{(N+1)}^m$,  ${\p}^{\B}={\p}_{{\xi}_k}$, 
then $b_{\B}\in S_{(N)}^m\Rightarrow  
{\tilde b}_{\B}=h^{1/2}i{\p}^{\B}b_{\B}\in S_{(N-1)}^m$ due to Lemma 3.1(b).  
Thus (4.8) gives   
$t^nJ^h_t(b)=t^nJ^h_t(b_0)+t^{n-1}J^h_t({\tilde b})$ with $b_0\in S^m_{(N)}$, 
${\tilde b}\in S^m_{(N-1)}$ and we complete the proof 
 using the assertion a) with $b_0$, ${\tilde b}$ instead of $b$.      
 $\triangle$ 
\par \bigskip 
{\bf Proof of (2.19).} Writing the terms 
$t^{k-1}J^h_t(q^{\; \circ}_{{\bar N},k})$ with $k\in \{ 2,...,{\bar N}\}$ 
as described in Lemma 4.1(b) (for $N=n=k-1$ and $b=q^{\; \circ}_{{\bar N},k}
\in {\check S}_{(k)}^0$), 
we can express (4.6) in the form 
 $$
 {\rm tr\,}Q^h_{\bar N}(t){\tilde L}_h^*\; =J_t^h({\check l}_h)+
\sum_{1\le k\le {\bar N}}t^kJ^h_t(b_{{\bar N},k}) \eqno (4.9)
 $$  
with some $b_{{\bar N},k}\in S_{(0)}^0$. Changing the order of integrals we 
find    
 $$ 
 \int_{\R}{dt\over 2\pi h}\, f_h^Z(t)J^h_t({\check l}_h)\, =\int_{{\R}^{2d}}
 {dv\over (2\pi h)^d}\, {\check l}_h(v)\; {\tilde f}^Z_h({\hat p_h}(v)).  
 $$
It remains to consider the terms of (4.9) with $k\in \{ 1,...,{\bar N}\}$ 
and to estimate    
 $$ 
 \int_{\R}{dt\over 2\pi h}\, f^Z_h(t)t^kJ^h_t(b_{{\bar N},k})
 \; =\int_{{\R}^{2d}}{dv\over (2\pi h)^d}\; b_{{\bar N},k,h}(v)
 i^{-k}h^k({\tilde f}_h^Z)^{(k)}({\hat p_h}(v)). \eqno (4.10) 
 $$
However 
 $k\ge 1 \Rightarrow h^k({\tilde f}_h^Z)^{(k)}(\la )={\tilde \g}_1^{(k-1)}
 \! \! \left( {  \la -E_1\over h} \right) -
 {\tilde \g}_1^{(k-1)}\! \! \left( {  \la -E_2\over h} \right)$ 
 \T  and we have 
 $|{\tilde \g}_1^{(k-1)}(\la )|\le C_{k,N}(1+|\la |)^{-N}$ for every 
$\la \in \R$, $N\in \N$, hence (4.10) can be estimated by  
$$
 \sum_{1\le j \le 2} h^{-d}\int_{{\overline \G}_{\! 0}}dv\; C'_{k,N}
{ \left( {  1+ { |{\hat p}_h(v) -E_j|\over Ch } 
} \right)  }^{\! \! -N}, 
 $$
where we can choose $C\ge 1$ such that 
$\sup_{v\in {\overline \G}_{\! 0}}|a_0(v)-{\hat p}_h(v)|\; \le {1\over 2}Ch$.\par  
It is clear that the region 
$\{ v\in {\R}^{2d}: |{\hat p}_h(v)-E_j|\ge {1\over 3}h^{1-\e}\}$ 
gives a contribution $O(h^{N\e -d})$ 
and it remains to consider the regions 
 $$
 {\G}_{\! E_j}^{h,n}=\{ v\in {\R}^{2d}: C(n-1/2)h\le {\hat p}_h(v)-E_j\le 
 C(n+1/2)h\} , \eqno (4.11)
 $$ 
where $n\in \Z$ is such that $|n|<h^{-\e}/(2C)$. 
However     
 $$
 {\G}_{\! E_j}^{h,n}\subset {\widetilde \G}_{\! E_j}^{h,n}=\{ v\in {\R}^{2d}: 
 C(n-1)h\le a_0(v)-E_j\le C(n+1)h\} 
 $$ 
and to complete the proof we observe  that
 $$
 \sum_{|n|<h^{-\e}/(2C)} \int_{ {\G}_{\! E_j}^{h,n}} dv
 { \left( {  1+ { |{\hat p}_h(v) -E_j|\over Ch } } \right)  }^{\! \! -N}
 \! \! \le  \left( { 3\! +\! \sum_{n=1}^{\infty}{2\over n^N} } \right)   \,
 \sup_{|n|<h^{-\e}/(2C)} {\rm vol}\, {\widetilde \G}_{\! E_j}^{h,n}
 $$
can be estimated by $C_{\e}{\cal R}^{\e ,a_0}_{E_j}(h)$. $\triangle$ 
\par \bigskip 
 Our proof of (2.18) will use
\begin{prop} Let ${\de}_0+{1\over 2}<\mu <1$ and 
$0<\ka \le \min \{ \mu -{\de}_0-{1\over 2},\; {1-\mu\over 2}\}$. Then    
 $$
 \sup_{\{ t\in \R:\; |t|\le h^{1-\mu}\} } |\, {\rm tr\;}  
 (Q^h_{\bar N}(t)-{\check L}_h{\rm e}^{itP_h/h}){\tilde L}_h^*\; |\T \le \, 
 C_{\bar N}h^{{\bar N}\ka -5d-1}.  \eqno (4.12)
 $$
\end{prop}
To begin we describe the form of 
${\tilde Q}^h_{\bar N}(t)$. Since $P=P^*=p(x,hD)^*$, we have 
 $$ 
 Q^h_{\bar N}(t)P_h={\rm Op}_t^h 
 [ q_{{\bar N},{\bar N}}(t,x,\xi ){\overline {p(y,\xi )} } ] 
 $$ 
and  the standard Taylor's development of ${\overline {p(\cdot ,\xi )} }$ 
in $x$, followed by integrations by parts based on
 $(x-y)^{\a}{\rm e}^{i(x-y)\xi /h}=
(-ih)^{|\a |}{\p}^{\a}_{\xi}({\rm e}^{i(x-y)\xi /h}),$ 
gives  
  $$
 {\tilde Q}^h_{\bar N}(t)={\rm Op}^h_t
 [ ({\tilde {\cal P}_{\bar N}}q_{{\bar N},{\bar N}})(t,x,\xi )+
r_{\bar N}(t,x,\xi ,y)]  \eqno (4.13)
  $$
  with 
  $$
  r_{{\bar N},h}(t,x, \xi ,y)=h^{-1}
{\rm e}^{-it{\hat p}_h(x, \xi )/h}({\bar N}+1)\int_0^1 d\si \; 
(1-\si )^{\bar N} {\tilde r}_{{\bar N},\si ,h}(t,x,\xi ,y),  
  $$
  $${\tilde r}_{{\bar N},\sigma ,h}(t,x,\xi ,y)=$$
  $$
  \sum_{|\a |={\bar N}+1} {(-ih)^{|\a |}\over \a !}  {\p}_{\xi}^{\a}  
\!{ \Bigl( (q_{{\bar N},{\bar N},h}(t){\rm e}^{it{\hat p}_h/h})(x,\xi )\,
 {\overline  { {\p}^{\a}_x p_h(x+\si (y-x), \xi )} } \Bigr) } .\eqno (4.14) 
  $$ 
We describe the properties of $r_{\bar N}(t)$ introducing new classes of 
symbols. \par 
We define ${\tilde S}^m_{\mu ,(N)}$ for $N\in {\N}$ 
setting ${\tilde S}^m_{\mu ,\; (0)}={\tilde S}^m_{(0)}$ and writing 
$b\in {\tilde S}^m_{\mu ,(N+1)}$ if and only if 
 $$
 b_h(x,\xi ,y)=b_{0,h}(x,\xi ,y)+\sum_{|\B |=1}b_{\B ,h}(x,\xi ,y)
 h^{-\mu /2}{\p}^{\B}{\hat p}_h(x,\xi ) \eqno (4.15)
 $$ 
holds with some $b_{\B}\in {\tilde S}^m_{\mu ,\; (N)}$ for 
$\B \in {\N}^{2d}$ satisfying $|\B |\le 1$. Then 
 $$
 {\tilde S}^m_{1,\; (N)}\subset {\tilde S}^{m+(1-\mu )N/2}_{\mu ,\; (N)} 
 \eqno (4.16)
 $$
and we adopt the following convention: 
every $(b_h)_{h\in ]0;\; h_0]}\in S^m_{(N)}$ is identified with  
$(x,\xi ,y)\to b_h(x,\xi )$ defining an element of ${\tilde S}^m_{1 ,\; (N)}$.  
 \par \bigskip 
\T
 To deduce  (4.12) using (3.3) we consider \T 
 $t=h^{1-\mu}{\tilde t}$, \T  $\tau =h^{1-\mu}{\tilde \tau}$ \T  
and \T ${\cal V}=\{ ({\tilde t},{\tilde \tau})\in ([-t_0;\; t_0]\setminus \{ 0\})^2:
 \; 0\le {\tilde \tau}/{\tilde t}\le 1 \}$. We denote 
$U^h_s={\rm e}^{isP/h^{\mu}}$ and we observe that (4.12) follows from   
 $$
 \sup_{({\tilde t},{\tilde \tau})\in {\cal V}} \left| {\, {\rm tr\;}  
 {\left( { {\tilde Q}^h_{\bar N}(t) U^h_{\tilde \tau}{\tilde L}_h^*} 
 \right) }_{t=h^{1-\mu}({\tilde t}-{\tilde \tau}) }\; }\right |
 \T \le \, C_{\bar N}h^{{\bar N}\kappa -5d-1}.  \eqno (4.17)
 $$
However using the form of $q_{{\bar N},{\bar N}}(t)$ in (4.14) we obtain the 
expressions of the form considered in the proof of Lemma 3.2 and applying 
(3.19), (3.20) with $|\a |={\bar N}+1$ we find 
$r_{{\bar N}}(t)=\sum_{n=0}^{2{\bar N}}t^nr_{{\bar N},n}^{\circ}$ with
 $$
 r_{{\bar N},n}^{\circ}\in 
{\tilde S}_{1,(2{\bar N}+1)}^{-({\bar N}+1)(1/2-{\de}_0)}
 \subset {\tilde S}_{\mu ,\; (2{\bar N}+1)}^{-({\bar N}+1)(\mu -1/2-{\de}_0)} 
 \eqno (4.18)
 $$
[the inclusion follows from (4.16)].  
Due to $(3.23({\bar N}))$ we have
 $$
 t=h^{1-\mu}{\tilde t}\T \Longrightarrow \T ({\tilde {\cal P}_{\bar N}}
 q_{{\bar N},{\bar N}})(t)
 =\sum_{{\bar N}\le n\le 2{\bar N}}{\tilde t}^{\; n}\; 
 {\tilde q}^{\; (\mu )}_{{\bar N},{\bar N},n} \eqno (4.19)
 $$
with 
${\tilde q}^{\; (\mu )}_{{\bar N},{\bar N},n,h}=h^{n(1-\mu )}
{\tilde q}^{\; \circ}_{{\bar N},{\bar N},n,h}$ belonging to 
${\tilde S}^{-n(1-\mu )}_{1,(n+1)}\subset 
{\tilde S}^{-n(1-\mu )/2}_{\mu ,\; (n+1)}$.\par  
Combining (4.18) and (4.19) we find the expression 
 $$
 t=h^{1-\mu}{\tilde t}\T \Longrightarrow \T {{\tilde Q}^h_{\bar N}(t)}
 =\sum_{0\le n\le 2{\bar N}}{\tilde t}^{\; n}\,  {\rm Op}^h_t[\; (1-
 {\tilde \tau}/{\tilde t}\; )^n\; b_{{\bar N},n}^{\, (\mu )}\; ], \eqno (4.20) 
 $$
with $b_{{\bar N},n}^{\, (\mu )}\in 
{\tilde S}^{-{\bar N}\ka }_{\mu ,\; (n+1)}$ similarly as in 
the formula (4.14) of [26] and 
 following [26] we denote 
  $$
  J^h_{{\tilde t},{\tilde \tau}}(b,Y)={\rm tr}\,  \Bigl( 
  {\rm Op}^h_t[b]U^h_{\tilde \tau}Y^h_{{\tilde t},{\tilde \tau}} 
  {\Bigr)}_{t=h^{1-\mu}({\tilde t}-{\tilde \tau})} \eqno (4.21)
  $$
if $b\in {\tilde S}^m_{\mu ,(N)}$ and 
$Y=(Y^h_{{\tilde t},{\tilde \tau}})_{(h,{\tilde t},{\tilde \tau}) 
\in ]0;\; h_0]\times {\cal V}}\subset B(L^2({\R}^d))$. 
Due to (4.20) the estimate (4.17) follows from  
 $$
 \sup_{({\tilde t},{\tilde \tau})\in {\cal V}} |\, {\tilde t}^{\; n}
J^h_{{\tilde t},{\tilde \tau}}(b_{{\bar N},n}^{\, (\mu )},
 {\tilde L}^*) |\T \le \; C_{\bar N}h^{{\bar N}\kappa -5d-1}.   \eqno (4.22)
 $$
Next we observe that the properties of operators $P_h$ given in Lemma 2.1 
allow us to follow the reasoning of Sections 5-6 in [26]. 
More precisely: let $T_j=h^{\mu /2-1}x_j$, 
$T_{-j}=h^{\mu /2}{\p}_{x_j}$, $P_{\pm j}=[ih^{-\mu}P,\; T_{\pm j}]$ 
for $j\in \{ 1,...,d\}$ and 
write $B\in {\Psi}^m_{{\de}_0}$ if and only if 
$B=(b_h(x,hD))_{h\in ]0;\; h_0]}$ holds with  
$(b_h)_{h\in ]0;\; h_0]}\in S^m_{{\de}_0}$. 
Using $\mu >{\de}_0+{1\over 2}>2{\de}_0$ it is easy to check  
that for every $B\in {\Psi}^0_{{\de}_0}$ 
we have $[B,T_{\pm j}]\in {\Psi}^0_{{\de}_0}$ and 
$[B,P_{\pm j}]\in {\Psi}^{-\ka}_{{\de}_0}$ with $\ka >0$. 
Thus it is easy to check that using $T_{\pm j}$, $P_{\pm j}$ 
as above and ${\Psi}^0_{{\de}_0}$ instead of ${\Psi}^0$ in the definition 
of ${\cal Y}$, we can follow the reasoning of
 the proof of Proposition 4.2 of [26] and  we obtain 
  
\begin{prop} Let $K=K({\bar N},n)\in \N$ be large enough. Then 
 one can write  
  $$
 {\tilde t}^{\; n}\; J^h_{{\tilde t},{\tilde \tau}}(b_{{\bar N},n}^{\, (\mu )}, 
 {\tilde L}^*)=\sum_{1\le k\le K} 
J^h_{{\tilde t},{\tilde \tau}}(b_{k,n},Y_{k,n})\, 
\T {\it for}\T ({\tilde t},{\tilde \tau})\in {\cal V},  \eqno (4.23)
  $$  
where $b_{k,n}\in {\tilde S}^{-{\bar N}\ka }_{(1)}$ (for $k\in \{ 1,...,K\}$) 
and $Y_{k,n}=(Y^h_{{\tilde t},{\tilde \tau},n,k}
)_{(h,{\tilde t},{\tilde \tau})\in ]0;\; h_0]\times {\cal V}}$ 
is a bounded subset of $B(L^2({\R}^d))$ (for $k\in \{ 1,...,K\}$). 
  \end{prop}
{\bf Proof of Proposition 4.2.} We observe that (4.22) follows from 
Proposition 4.3 similarly as in [26]. Indeed, using the 
expression (4.23) it suffices to write  
 $$
 |J^h_{{\tilde t},{\tilde \tau}}(b_{k,n},Y_{k,n})|\le 
||{\rm Op}^h_t[b_{k,n}]_{t=h^{1-\mu}({\tilde t}-{\tilde \tau})}\, ||_{\rm tr}\; 
||Y^h_{{\tilde t},{\tilde \tau},n,k}||
 \le Ch^{{\bar N}\ka -5d-1/2}, 
 $$
where we used (4.3) with $q=b_{k,n}\in {\tilde S}^{-{\bar N}\ka }_{(1)}
\subset {\tilde S}^{1/2-{\bar N}\ka }_{(0)}$. 
$\triangle$ 
 \par \bigskip 
To complete the proof of Proposition 2.3 it remains to use 
\begin{prop} If $\mu >{\de}_0+{1\over 2}$ then  
 $$ 
 \sup_{\{ t\in \R :\; h^{1-\mu}\le |t|\le t_0\} } |\, {\rm tr\,}{\check L}_h\; 
{\rm e}^{itP_h/h}{\tilde L}_h^*\, |\T =\, O(h^{\infty}),\eqno (4.24)
 $$
 $$ 
 \sup_{\{ t\in \R :\; h^{1-\mu}\le |t|\le t_0\} } |\, {\rm tr\,} 
Q^h_{\bar N}(t){\tilde L}_h^*\, |\T =\, O(h^{\infty}). \eqno (4.25)
 $$
\end{prop}

\section{Proof of Proposition 4.4} 

{\bf Proof of (4.25).} Let $\kappa$ be as in Proposition 4.2 and 
 $$
 b=(b_h)_{h\in ]0;\; h_0]}\in S^m_{{\de}_0}({\R}^{2d}) 
\T {\rm with}\T {\rm supp\;}b_h\subset {\check \G}({\bar c}h^{{\de}_0})
 \cap {\overline \G}_{\! 0}. \eqno (5.1)
 $$  
We will show that 
 $$
 \sup_{\{ t\in \R :\; h^{1-\mu}\le |t|\le t_0\} }  
 |J^h_t(b)|\, =\, O(h^{n\ka -m}) \eqno (5.2(n)) 
 $$
holds for every $n\in \N$, which ensures (4.25) due to (4.6). \par 
Reasoning by induction we assume that the assertion holds for a given $n\in \N$. 
Using the cut-off functions from the proof of Lemma 2.4, 
it is easy to see that the assumptions (5.1) ensure the existence of 
$b_{\B}\in S^m_{{\de}_0}({\R}^{2d})$ such that 
$b=\sum_{|\B |=1}b_{\B}h^{-{\de}_0}{\p}^{\B}{\hat p}$ 
and ${\rm supp\;}b_{\B}\subset {\rm supp\;}b$.  
The integration by parts gives 
 $$
 J^h_t(b)=\sum_{|\B |=1}t^{-1}J^h_t(h^{1-{\de}_0}i{\p}^{\B}b_{\B}). 
\eqno (5.3)
 $$ 
The induction hypothesis applied to 
$h^{1-{\de}_0}{\p}^{\B}b_{\B}\in S^{m+2{\de}_0-1}_{{\de}_0}({\R}^{2d})$ 
ensures 
 $$
 h^{1-\mu}\le |t|\le t_0\; \Rightarrow |t^{-1}J^h_t(h^{1-{\de}_0}
{\p}^{\B}b_{\B})| \le Ch^{-(1-\mu )+n\ka -m-2{\de}_0+1}. \eqno (5.4)
 $$  
Since  $\mu -2{\de}_0>\mu -{\de}_0-{1\over 2}\ge \ka$, it is clear that 
(5.3--4) imply (5.2($n+1$)). $\triangle$ 
 \par 
 \bigskip 
 Let ${\vt}^h_t\! :{\R}^{2d}\to {\R}^{2d}$ be the Hamiltonian 
flow of ${\hat p}_h$, i.e. $t\to {\vt}^h_t(v)$ satisfies
 $$
 \hbox{${d\over dt}$} {\vt}^h_t(v)={\cal J}\nabla {\hat p}_h({\vt}^h_t(v)), 
 \hspace{1cm} {\vt}^h_t(v)|_{t=0}\, =v, 
 $$    
where ${\cal J}= {\pmatrix { 0_{{\R}^d} & I_{{\R}^d}\cr 
-I_{{\R}^d} & 0_{{\R}^d}\cr} }$. \T

\begin{lem} Assume that $t_0>0$ is small enough. Then 
 $$ 
 -t_0\le t\le t_0\Rightarrow 
|{\vt}^h_t(v)-v|\ge |t\nabla {\hat p}_h(v)|/2. 
  \eqno (5.5)
 $$
\end{lem} 
{\bf Proof.} Set 
  $M^h_t(v)=\int_0^1ds\; {\cal J}\nabla d{\hat p}_h(v+s({\vt}^h_t(v)-v))$.  
Then 
 $$
 {\cal J}\nabla {\hat p}_h({\vt}^h_t(v))-{\cal J}\nabla {\hat p}_h(v)=
 M^h_t(v)({\vt}^h_t(v)-v) \eqno (5.6)
 $$
 and $\hbox{${d\over dt}$}({\vt}^h_t(v)-v)=M^h_t(v)({\vt}^h_t(v)-v)+
 {\cal J}\nabla {\hat p}_h(v).$ 
Therefore introducing the solution of the linear homogeneous system
 $$
  \hbox{${d\over dt}$} R^h_{t,\tau}(v)=M^h_t(v)R^h_{t,\tau}(v) \hspace{1cm} 
 R^h_{t,\tau}(v)|_{t=\tau}\, =I 
 $$
we obtain 
${\vt}^h_t(v)-v=\int_0^td\tau R^h_{t,\tau}(v){\cal J}\nabla {\hat p}_h(v)$, 
which ensures 
 $$
 -1\le t\le 1\Rightarrow 
 |{\vt}^h_t(v)-v|\le C_1|t\nabla {\hat p}_h(v)|. \eqno (5.7)
 $$   
Using \T $\int_0^t d\tau \; {d\over d\tau} ({\vt}^h_{\tau}(v)-v-\tau {\cal J}\n 
{\hat p}_h(v))=\int_0^t d\tau \; {\cal J}(\n{\hat p}_h({\vt}^h_{\tau}(v))
 -\n {\hat p}_h(v))$ \T 
and (5.6) we obtain 
 $$
 |{\vt}^h_t(v)-v-t{\cal J}\nabla {\hat p}_h(v)|\le \int_0^t d\tau \, 
|M^h_{\tau}(v)||{\vt}^h_{\tau}(v)-v|. \eqno (5.8) 
 $$
Using (5.7) to estimate the right hand side of (5.8) we obtain 
 $$
 -1\le t\le 1\Rightarrow |{\vt}^h_t(v)-v-t{\cal J}\nabla {\hat p}_h(v)|\le 
C_2t^2|{\cal J}\nabla {\hat p}_h(v)|.  \eqno (5.9)
 $$
Writing  
$ |{\vt}^h_{\tau}(v)-v|\ge |t{\cal J}\nabla {\hat p}_h(v)|-
 |{\vt}^h_t(v)-v-t{\cal J}\nabla {\hat p}_h(v)|$ and 
using (5.9) we obtain (5.5) if $|t|\le \min \{ 1,\; 1/(2C_2)\}$.
 $\triangle$ 
\par \bigskip 
\T For $h\in ]0;\; h_0]$ let ${\Omega}_h$ be a set of parameters. 
We say that the family 
$(b_{\omega ,h})_{(\omega ,h)\in {\Omega}_h\times ]0;\; h_0]}$ 
is bounded in $S^m_{\de}(\G )$ if and only if the estimates 
 $$
\sup_{\omega \in {\Omega}_h} \sup_{v\in {\G}_h} |{\p}^{\a}b_{\omega ,h}(v)|
 \le C_{\a }h^{-m-|\a |{\de}}  \eqno (5.10)
 $$
hold for all $\a \in {\N}^{2d}$.

\begin{lem} Assume that $l,{\tilde l}\in S^0_{\de}({\R}^{2d})$ satisfy 
 $$
 {\rm dist}(\; {\rm supp\,} {\tilde l}_h, {\rm supp\,} (1-l_h)\, )\, \ge 
 h^{\de}. \eqno (5.11)   
 $$
Let ${\tilde L}_h={\tilde l}_h(x,hD)$ and  
$L_h(t)=(l_h\! \circ \! {\vt}^h_t)(x,hD)$. Then 
 $$
 ||(I-L_h(t)){\rm e}^{itP_h/h}{\tilde L}_h||\, =\, O(h^{\infty}).\eqno (5.12)
 $$
\end{lem} 
{\bf Proof.} We observe that for $\a \in {\N}^{2d}$ such that $|\a |\le 1$, 
the matrix elements of 
$({\p}^{\a}{\vt}^h_t)_{(t,h)\in [-t_0;\; t_0]\times ]0;\; h_0]}$ 
are bounded families in $S^0_{{\de}_0}({\overline \G}_{\! 0})$ 
and it is easy to check that for every   
$l\in S^m_{{\de}_0}({\R}^{2d})$, the family  
$(l_h\! \circ \! {\vt}^h_t)_{(t,h)\in [-t_0;\; t_0]\times ]0;\; h_0]}$ 
is bounded  in $S^m_{{\de}_0}({\overline \G}_{\! 0})$. 
Then the properties of operators $P_h$ 
given in Lemma 2.1 allow us to follow the proof of Lemma 5.1 and 
Proposition 5.2 of [24] with $|t|\le t_0$ instead of 
$|t|\le 2h^{{\de}_0}$. $\triangle$ 
\par \bigskip     
{\bf Proof of (4.24).} We can consider a family of 
balls $\{ B({\bar v}_{n,h},h^{\de})\} {}_{n\in \{ 1,...,N(h)\} }$ covering 
${\rm supp\;}{\check l}\subset {\check \G}({\bar c}h^{{\de}_0})\cap 
{\overline \G}_{\! 0}$ with ${\bar v}_{n,h}\in {\check \G}({\bar c}h^{{\de}_0})\cap 
{\overline \G}_{\! 0}$ for $n\in \{ 1,...,N(h)\}$ 
and $N(h)\le Ch^{-2d\de}$. 
Using a suitable partition of unity 
we  decompose 
 $$
 {\check l}_h=\sum_{1\le n\le N(h)}{\check l}_{n,h}\T \T {\rm with}\T 
 {\rm supp\;}{\check l}_{n,h}\subset B({\bar v}_{n,h},h^{\de}), 
 $$
where  $({\check l}_{n,h})_{(n,h)\in \{ 1,...,N(h)\} \times ]0;\; h_0]}$ 
is bounded in $S^0_{\de}({\R}^{2d})$. \par  

  Let $({\tilde l}_{n,h})_{(n,h)\in \{ 1,...,N(h)\} \times ]0;\; h_0]}$, 
 $(l_{n,h})_{(n,h)\in \{ 1,...,N(h)\} \times ]0;\; h_0]}$  
be  bounded in $S^0_{\de}({\R}^{2d})$ and such that
 ${\rm supp\;}{\tilde l}_{n,h}\subset B({\bar v}_{n,h},2h^{\de}),\T \T 
 {\rm supp\;}l_{n,h}\subset B({\bar v}_{n,h},4h^{\de})$. 
Using the trace cyclicity and assuming ${\tilde l}_{n,h}=1$ on 
$B({\bar v}_{n,h},h^{\de})$ we find 
 $$
 {\rm tr\;}{\check L}_h{\rm e}^{itP_h/h}{\tilde L}^*_h=\sum_{1\le n\le N(h)} 
 {\rm tr\;}{\rm e}^{itP_h/h}{\tilde L}^*_h{\check L}_{n,h}
 $$
 $$
 =\sum_{1\le n\le N(h)}{\rm tr\;}{\rm e}^{itP_h/h}{\tilde L}_{n,h}
 {\tilde L}^*_h{\check L}_{n,h}+O(h^{\infty}), \eqno (5.13)
 $$
where we have denoted ${\check L}_{n,h}={\check l}_{n,h}(x,hD)$ and 
 ${\tilde L}_{n,h}={\tilde l}_{n,h}(x,hD)$. \par 
We assume  $l_{n,h}=1$ on 
$B({\bar v}_{n,h},3h^{\de})$ and introduce 
$L_{n,h}(t)=(l_{n,h}\! \circ \! {\vt}^h_t)(x,hD)$. Since the 
assertion of Lemma 5.2 still holds with $l_{n,h}$ and ${\tilde l}_{n,h}$ 
instead of $l_h$ and ${\tilde l}_h$,  (5.13) can be written as 
 $$
 \sum_{1\le n\le N(h)}{\rm tr\;}L_{n,h}(t){\rm e}^{itP_h/h}{\tilde L}_{n,h}
 {\tilde L}^*_h{\check L}_{n,h}+O(h^{\infty}). \eqno (5.14)
 $$
To complete the proof it suffices to check that 
 $$
 h^{1-\mu}\le |t|\le t_0\Rightarrow {\rm supp\;}{\check l}_{n,h}\cap 
 {\rm supp\;}(l_{n,h}\! \circ \! {\vt}^h_t)\, =\emptyset \eqno (5.15)
 $$  
holds for a certain $\de \in [0;\; 1/2[$. Indeed, (5.15) ensures 
$||{\check L}_{n,h}L_{n,h}(t)||=O(h^{\infty})$ and (5.14) is $O(h^{\infty})$ 
due to the trace cyclicity. \par 
\T To obtain (5.15) we observe that 
 $$
 {\rm supp\;}(l_{n,h}\! \circ \! {\vt}^h_t)\subset {\vt}^h_{-t}(B({\bar v}_{n,h},
4h^{\de}))\subset B({\vt}^h_{-t}({\bar v}_{n,h}),C_0h^{\de})\Rightarrow 
 $$
 $$
 {\rm dist}(\; {\rm supp\,} {\check l}_{n,h}, {\rm supp\;}(l_{n,h}\! \circ \! 
{\vt}^h_t)\, )\ge |{\vt}^h_{-t}({\bar v}_{n,h})-{\bar v}_{n,h}|
-(1+C_0) h^{\de}. \eqno (5.16)   
 $$
\par \smallskip 
Since ${\bar v}_{n,h}\in {\check \G}({\bar c}h^{{\de}_0})\cap 
{\overline \G}_{\! 0}\Rightarrow 
|\nabla {\hat p}_h({\bar v}_{n,h})| \ge {\bar c}h^{{\de}_0}/2$,  
Lemma 5.1 ensures 
 $$
h^{1-\mu}\le |t|\le t_0\Rightarrow |{\vt}^h_{-t}({\bar v}_{n,h})-{\bar v}_{n,h}|
\ge  |t\nabla {\hat p}_h({\bar v}_{n,h})|/2\ge {\bar c}
 h^{1-\mu +{\de}_0}/4 
 $$ 
and (5.16) implies (5.15) if we take $\de \in ]1-\mu +{\de}_0;\; 1/2[$. 
$\triangle$

\section{End of the proof of Theorem 1.2}
 Since $d^2a_0$ is continuous, 
we can assume that $c>0$ is small enough to ensure 
 $$
 v\in {\cal C}^{a_0}_E(c)\, \Rightarrow \, {\rm rank}(d^2a_0(v))\ge 2, 
 \eqno (6.1)
 $$
where 
$$
{\cal C}^{a_0}_E(c)=\{ v\in {\R}^{2d}:\, |a_0(v)-E|+ 
|\n a_0(v)|\le 2c \, \} .
$$ 
 Until now we have proved 
$$
 \left| { {\rm tr\;}{\tilde f}^Z_h(P_h){\check L}_h\,-
 \int_{{\R}^{2d}} {dv \over (2\pi h)^d}\, {\tilde f}^Z_h
 ({\hat p}_h(v )){\check l}_h(v) } \right|
  \le C_{\e}\, h^{-d} \sum_{1\le j\le 2}\! {\cal R}^{\e ,a_0}_{E_j}(h)
   \eqno (6.2)
 $$
and we can deduce Theorem 1.2 from Lemma 2.2 if we show 

\begin{prop} Let $l$ be as in Lemma $2.2$,  ${\hat l}$ as in Lemma $2.4$ and 
 ${\hat L}_h={\hat l}_h(x,hD)$.  
 If $(6.1)$ holds and   
 $$
 \hbox{${1\over 2}$} \left( {1-\hbox{${1\over 4m_0-1}$} } \right) 
 <{\de}_0 < \hbox{${1\over 2}$},  \eqno (6.3)
 $$ 
 then 
 $$
 {\rm tr\;}{\hat L}_h{\tilde f}^Z_h(P_h)\, =\; o(h^{1-d}), \eqno (6.4)
 $$ 
 $$
 \int_{{\R}^{2d}} {\hat l}_h({\tilde f}^Z_h\! \circ \! {\hat p}_h)\, =\; o(h).
 \eqno (6.5)
 $$ 
\end{prop} 

The proof of Proposition 6.1 uses the following trace norm estimate
 
\begin{lem} Let ${\hat L}_h={\hat l}_h(x,hD)$ with 
${\hat l}\in S^0_{{\de}_0}({\R}^{2d})$ and denote 
 $$
 {\G}_{\! h}={\rm supp\,}{\hat l}_h\, +B(0,h^{{\de}_0})= \{ v\in {\R}^{2d}:\, 
{\rm dist}(v,{\rm supp\,}{\hat l}_h)\, < h^{{\de}_0}\, \} .  \eqno (6.6)
 $$ 
Then \T $||{\hat L}_h||_{\rm tr}\le Ch^{-d\, }{\rm vol}\, {\G}_h.$
\end{lem}
{\bf Proof.} Let $B_h=b_h(x,hD)$ with $b\in S^0_{{\de}_0}({\R}^{2d})$. 
It is well known (cf. e.g. [21]) that the Hilbert-Schmidt norm 
 $$
 ||B_h||_{\rm HS}=(2\pi )^{-d} 
 {\left( { \int_{{\R}^{2d}}|b_h(v)|^2h^{-d}dv} \right)}^{1/2} 
\le Ch^{-d/2}{ \left( { {\rm vol}\, [{\rm supp\,}b_h] } \right) }^{1/2}.
 $$ 
We can assume $b_h=1$ on ${\rm supp\,}{\hat l}_h$ and 
${\rm supp\,}b_h\subset {\G}_h$. Therefore we have \par 
 $||{\hat L}_h(I-B_h)||_{\rm tr}=O(h^{\infty})$
  and we complete the proof writing 
 $$
 ||{\hat L}_hB_h||_{\rm tr}\le ||{\hat L}_h||_{\rm HS}||B_h||_{\rm HS}\le 
 Ch^{-d}{ \left( { {\rm vol}\, [{\rm supp\,}{\hat l}_h]\cdot {\rm vol}\, 
[{\rm supp\,}b_h] } \right) }^{1/2}. \eqno \triangle
 $$
Instead of Proposition 6.1 it suffices to show

\begin{prop} Let $l$ be as in Lemma $2.2$ and 
 ${\hat l}$ as in Lemma $2.4$. If ${\G}_{\! h}$ is given by $(6.6)$, then 
 $$
 {\rm vol}\, {\G}_{\! h} \, =O(h^{{\de}_0(4m_0-1)/(2m_0-1)}). \eqno (6.7)
 $$
\end{prop} 
Indeed, since (6.3) ensures ${\de}_0(4m_0-1)/(2m_0-1)>1$, it is 
obvious that (6.7) implies (6.5) and using Lemma 6.2 we obtain similarly (6.4).  
\par \bigskip 
To begin the proof of Proposition 6.3 we introduce the following notation: 
 $$
 \cases { {\p}_j={\p}_{x_j}\TT {\rm if}\T j\in J_+=\{ 1,...,d\} \cr  
 {\p}_j={\p}_{{\xi}_{-j}}\TT {\rm if}\T j\in J_-=\{ -1,...,-d\} } .     
 $$
Let ${\bar v}\in {\R}^{2d}$ be such that ${\rm rank}(d^2a_0({\bar v}))\ge 2$. 
Then there exist 
$j(1,{\bar v}),j(2,{\bar v})\in J=J_+\cup J_-$  such that 
 $\nabla {\p}_{j(k,{\bar v})}a_0({\bar v})\ne 0$ for $k\in \{ 1,2\}$ and the angle 
 $\an (\n {\p}_{j(1,{\bar v})}a_0({\bar v}),\n {\p}_{j(2,{\bar v})}a_0({\bar v}))
 \ne 0.$   
Therefore we can find two linearly independent vectors 
 $e_{1,{\bar v}},e_{2,{\bar v}}\in {\R}^{2d}$ satisfying 
 $$
 e_{k,{\bar v}}\cdot \n {\p}_{j(k,{\bar v})}a_0({\bar v})\, 
 ={\theta}_{k,{\bar v}}>0, 
\eqno (6.8)
 $$ 
 $$
 e_{k,{\bar v}}\notin \{ (0,\xi )\in {\R}^{2d}:\, \xi \in {\R}^d\} 
\eqno (6.9)
 $$
for $k\in \{ 1,2\}$. Since $d\ge 2$ we can find 
 $$
 e_{3,{\bar v}}\in \{ (0,\xi )\in {\R}^{2d}:\, \xi \in {\R}^d \}  \eqno (6.10)
 $$
such that the system $(e_{k,{\bar v}})_{k\in \{ 1,2,3\} }$ 
is linearly independent. 

\begin{lem} Let ${\bar v}\in {\R}^{2d}$ 
be such that ${\rm rank}(d^2a_0({\bar v}))\ge 2$  and let 
${\e}_{\bar v}>0$ be small enough. 
If ${\hat \G}_{\! h}$ is given by $(2.20)$ and $v\in {\R}^{2d}$, 
then the Lebesgue measure of 
 $$
 \{ s\in \R :\, se_{k,{\bar v}}+v\in B({\bar v},{\e}_{\bar v})\cap 
{\hat \G}_{\! h}\,  \} \eqno (6.11(k))
 $$ 
can be estimated by $C_{\bar v}h^{{\rho}_k}$, where 
  $$
{\rho}_1={\rho}_2={\de}_0, \T \T \T {\rho}_3={\de}_0/(2m_0-1)  \eqno (6.12)
 $$ 
and the constant $C_{\bar v}$ is independent of 
$(h,v)\! \in ]0;\; h_0]\times {\R}^{2d}$.  
\end{lem} 

Let us check that Proposition 6.3 follows from Lemma 6.4. 
To begin we observe that  (2.10) gives 
 $$
 v\in {\G}_{\! h}\subset {\rm supp\,}l_h\,  +B(0,h^{{\de}_0}) \Rightarrow 
 |a_0(v)-E|\le c+Ch^{{\de}_0} 
 $$ 
 and  using ${\G}_{\! h}\subset {\rm supp\,}{\hat l}_h\subset {\hat \G}_{\! h}$  
 we can choose $h_0>0$ small enough to ensure 
 ${\G}_{\! h}\subset {\cal C}^{a_0}_E(c)$ for $h\in ]0;\; h_0]$. 
The family $\{ B({\bar v},{\e}_{\bar v})\} {}_{{\bar v}\in 
{\cal C}^{a_0}_E(c) }$ contains a finite covering of ${\cal C}^{a_0}_E(c)$ 
and it suffices to show that for every ${\bar v}\in {\cal C}^{a_0}_E(c)$ 
one has 
 $$
 {\rm vol}\, { B({\bar v},{\e}_{\bar v})\cap {\hat \G}_{\! h} } \,  
\le {\bar C}_{\bar v}h^{{\rho}_1+{\rho}_2+{\rho}_3} . \eqno (6.13)
 $$ 
The system $(e_{k,{\bar v}})_{k\in \{ 1,2,3\} }$ 
can be completed to a basis 
$(e_{k,{\bar v}})_{k\in \{ 1,\dots ,2d\} }$ and 
let $(e^{\, *}_{k,{\bar v}})_{k\in \{ 1,\dots ,2d\} }$ 
denote the dual basis in $({\R}^{2d})^*$. 
Let $(e_k)_{k\in \{ 1,\dots ,2d\} }$ denote the 
canonical basis of ${\R}^{2d}$ and 
let $(e^{\, *}_k)_{k\in \{ 1,\dots ,2d\} }$ be its dual. 
 Then the set (6.11($k$)) has the form
 $$
 e^{\, *}_{k,{\bar v}}(B({\bar v},{\e}_{\bar v})\cap 
{\hat \G}_{\! h}\, -v)=e^{\, *}_k(W_{\bar v} (B({\bar v},{\e}_{\bar v})\cap 
{\hat \G}_{\! h}\, -v)),
 $$
where 
$W_{\bar v}\in {\rm Hom}({\R}^{2d},{\R}^{2d})$ is the matrix of the 
corresponding change of variables. The assertion of Lemma 6.4 allows us 
to estimate the measure of 
 $$
 \{ s\in \R :\, se_k+v\in W_{\bar v}( B({\bar v},{\e}_{\bar v})\cap 
{\hat \G}_{\! h})\, \}  \eqno (6.14)
 $$ 
by $C_{\bar v}h^{{\rho}_k}$ for $(h,v)\! \in ]0;\; h_0]\times {\R}^{2d}$, hence   
the Fubini's theorem gives
 $$
 {\rm vol}\, W_{\bar v}( B({\bar v},{\e}_{\bar v})\cap {\hat \G}_{\! h})
 \, \le C^3_{\bar v}h^{{\rho}_1+{\rho}_2+{\rho}_3}. \eqno (6.15) 
 $$
However 
${\rm vol}\, B({\bar v}, {\e}_{\bar v})\cap {\hat \G}_{\! h}\,  =
|\det W_{\bar v}|^{-1}{\rm vol}\, W_{\bar v}(B({\bar v},{\e}_{\bar v})\cap 
{\hat \G}_{\! h})$ and it is clear that  
 (6.13) follows from (6.15). 
\par \bigskip 
{\bf Proof of Lemma 6.4.}  
To begin we observe that the set (6.11(k)) is included in the interval
 $$
 {\D}_{k,{\bar v},v}=\{ s\in \R :\, se_{k,{\bar v}}+v\in 
B({\bar v},{\e}_{\bar v}) \,  \} .\eqno (6.16)
 $$ 
For $k\in \{ 1,\dots ,2d\}$ let $u_{k,{\bar v},v}:\, {\D}_{k,{\bar v},v}\to \R$ 
be defined by the formula
 $$
 u_{k,{\bar v},v}(s)={\p}_{j(k,{\bar v})}a_0(se_{k,{\bar v}}+v). \eqno (6.17)
 $$
By the definition of ${\hat \G}_{\! h}$ we find that 
the set (6.11($k$)) is included in 
 $$
 \{ s\in {\D}_{k,{\bar v},v}:\, -{\bar C}h^{{\de}_0}\le 
 u_{k,{\bar v},v}(s)\le {\bar C}h^{{\de}_0} \} .\eqno (6.18(k))
 $$ 
We claim that using (6.8) we can ensure 
 $$
s\in {\D}_{k,{\bar v},v}\Rightarrow \hbox{${d\over ds}$} u_{k,{\bar v},v}(s)
 >{\theta}_{k,{\bar v}}/2, \eqno (6.19)
 $$
for $k\in \{ 1,2\}$ if ${\e}_{\bar v}>0$ is fixed small enough. 
It is clear that (6.19) implies the fact that (6.18($k$)) defines an interval 
of length smaller than $2{\bar C}h^{{\de}_0}/{\theta}_{k,{\bar v}}$. \par   
In order to prove (6.19) we observe that
 $$
 \hbox{${d\over ds}$} u_{k,{\bar v},v}(s)=e_{k,{\bar v}}\cdot \! 
 \n {\p}_{j(k,{\bar v})}a_0(se_{k,{\bar v}}+v),
 $$
hence 
 ${d\over ds} u_{k,{\bar v},{\bar v}}(0)={\theta}_{k,{\bar v}}>0$. 
Moreover
 $$
 s\in {\D}_{k,{\bar v},v}\Rightarrow \left| {\hbox{${d\over ds}$} 
 u_{k,{\bar v},v}(s)-\hbox{${d\over ds}$} u_{k,{\bar v},{\bar v}}(0)} \right| 
 \le C|se_{k,{\bar v}}+v-{\bar v}|^{r_0}\le C{\e}_{\bar v}^{r_0}
 $$ 
and (6.19) follows if ${\e}_{\bar v}$ is such that 
$C{\e}_{\bar v}^{r_0}<{\theta}_{k,{\bar v}}/2$.
\par 
To estimate the measure of (6.11(3)) we observe that (6.18(3)) holds if we 
take $u_{3,{\bar v},v}(s)=e_{3,{\bar v}}\cdot \! \n a_0(se_{3,{\bar v}}+v)$,  
which is polynomial of degree $2m_0-1$ due to (6.10) and 
the ellipticity hypothesis (1.9) 
ensures $({d\over ds})^{2m_0-1}u_{3,{\bar v},v}(s)=
({d\over ds})^{2m_0-1}u_{3,{\bar v},v}(0)\ge c_0>0$. 
Thus to complete the proof of Lemma 6.4 it suffices to show 

 \begin{lem} Let $(F_{\o})_{\o \in \O}$ be a family of polynomials 
of order $m\in \N$ and 
 $$
 {\D}^{\! h}_{\o}=\{ s\in \R :\, |F_{\o}(s)|<(Ch)^{{\de}_0}\} .
 $$
If the $m$-th derivative satisfies $|F^{(m)}_{\o}(0)|\ge 1$, then 
the Lebesgue measure of ${\D}^{\! h}_{\o}$ can be estimated by 
$C_mh^{{\de}_0/m}$, where $C_m$ is independent of $\o \in \O$.   
\end{lem}
{\bf Proof.} We drop the index $\o$ and for $k\in \{ 0,...,m\}$ we set
 $$
 {\D}^{\! h}_k=\{ s\in \R :\, |F^{(k)} (s)|<(Ch)^{{\de}_k}\} ,
 $$
where ${\de}_k={\de}_0(1-k/m)$. We can write 
${\D}^{\! h}_0$ as the union of 
 $$
 {\D}^{\! h}_{(j_1,j_2,\dots ,j_m)}={\D}^{\! h}_0\cap {\D}^{\! h,j_1}_1\cap 
{\D}^{\! h,j_2}_2\cap \dots \cap {\D}^{\! h,j_m}_m,
 $$ 
where $j_k\in \{ 1,-1\}$, ${\D}^{\! h,1}_k={\D}^{\! h}_k$ and 
${\D}^{\! h,-1}_k=\R \setminus {\D}^{\! h}_k$. \par 
Since ${\de}_m=0$, the assumption  $|F^{(m)}(0)|\ge 1$ implies 
${\D}^{\! h}_m=\emptyset$ and for every $(j_1,...,j_m)\in {\{ 1,-1\}}^m$ 
we can find $k\in \{ 1,...,m\}$ such that 
 $$
 {\D}^{\! h}_{(j_1,j_2,\dots ,j_m)}\subset {\D}^{\! h}_{k-1}
 \setminus {\D}^{\! h}_k
 $$
However $F^{(k)}$ is a polynomial of order $m-k$, ${\D}^{\! h}_k$ 
is a union of at most $m-k$ intervals and 
 $$
 {\D}^{\! h}_{k-1}\setminus {\D}^{\! h}_k\, = \bigcup_{1\le j\le j(k,h)} 
 {\D}^{\! h}_{k,j}
 $$  
 where $j(k,h)\le 2(m-k)$ and ${\D}^{\! h}_{k,j}$ are intervals such that 
 $$
 s\in {\D}^{\! h}_{k,j}\Rightarrow |F^{(k-1)}(s)|<(Ch)^{{\de}_{k-1}} \T 
\hbox{and} \T |F^{(k)}(s)|\ge (Ch)^{{\de}_{k}}.
 $$
It is clear that the length of   ${\D}^{\! h}_{k,j}$ is 
less than $2(Ch)^{{\de}_{k-1}-{\de}_k}=2(Ch)^{{\de}_0/m}$. $\triangle$ 

 \section{Appendix} 
  
{\bf Proof of Lemma 2.1.} For $\a \in {\N}^d$ we denote 
${\g}^{\a}_h(x)=({\p}^{\a}\g )(h^{-{\de}_0}x)h^{-d{\de}_0}$. 
Further on $\a =\a ' +\a ''$ with $|\a ''|\le 2$ and dropping the indices 
$\nu ,{\bar \nu}$ we write 
 $$
 {\p}^{\a}a_h(x)={\p}^{\a '+\a ''}a_h(x)=\int_{{\R}^d} {\p}^{\a ''}{\! a}(y)
 {\g}_h^{ \a ' }(x-y)h^{-{\de}_0 |\a '|}\, dy. \eqno(7.1)
 $$
If $|\a |\le 2$ then it is clear 
that ${\p}^{\a}a_h=O(1)$ follows from (7.1) with $\a ' =0$. Further on 
we assume $|\a |\ge 3$ and $|\a ''|=2$. Then  
$|\a '|\ge 1\Rightarrow \int {\g}_h^{\a '}(x-y)dy=0$ and 
(7.1) still holds if ${\p}^{\a ''}{\! a}(y)$ is replaced by 
${\p}^{\a ''}{\! a}(y)-{\p}^{\a ''}{\! a}(x)$. Therefore 
 $$
 |{\p}^{\a}a_h(x)|\le \int_{{\R}^d}  |{\p}^{\a ''}{\! a}(y)-
 {\p}^{\a ''}{\! a}(x)|
 |{\g}_h^{\a '}(x-y)|h^{-{\de}_0|\a '|}\, dy \, \le 
 $$ 
 $$
 C \int_{{\R}^d}|y-x|^{r_0} |{\g}_h^{\a '}(x-y)|h^{-{\de}_0|\a '|}\, dy\, 
 =\, Ch^{(r_0-|\a '|){\de}_0}\, \int_{{\R}^d} |z|^{r_0} |{\g}^{(\a ')}(z)|\, 
 dz, \eqno (7.2)
 $$
i.e. we obtain ${\p}^{\a}a_h=O(h^{(r_0+2-|\a |){\de}_0})$ and (2.5) follows.  
 \par \smallskip 
If $|\a ''|=2$ then we can estimate
 $$
 {\p}^{\a ''}{\! a}_h(x)-{\p}^{\a ''}{\! a}(x)=\int_{{\R}^d}  
 ({\p}^{\a ''}{\! a}(y)-{\p}^{\a ''}{\! a}(x))
 {\g}^0_h(x-y)\, dy \eqno (7.3)
 $$ 
by $O(h^{r_0{\de}_0})$ similarly as the right hand side of (7.2) 
with $\a '=0$. \par 
In the next step we estimate (7.3) when $|\a ''|=1$. Since 
$\int y{\g}(y)\, dy\, =0$, we can replace 
${\p}^{\a ''}{\! a}(y)-{\p}^{\a ''}{\! a}(x)$ by  
  $$ 
 {\p}^{\a ''}{\! a}(y)-{\p}^{\a ''}{\! a}(x)-(y-x)\cdot \n {\p}^{\a ''}
 {\! a}(x) \eqno (7.4)
  $$ 
in the right hand side of (7.3). We can express (7.4) as 
 $$
 \int_0^1 \, ds \, 
 (y-x)\cdot (\n {\p}^{\a ''}{\! a}(x+s(y-x))- \n {\p}^{\a ''}{\! a}(x)), 
 \eqno(7.4')
 $$
 hence its absolute value is $O(|y-x|^{1+r_0})$ and  
 $$
 |{\p}^{\a ''}{\! a}_h(x)-{\p}^{\a ''}{\! a}(x)|\le C\int_{{\R}^d}  
 |y-x|^{1+r_0}|{\g}^0_h(x-y)|\, dy\, 
 =C_{r_0}h^{(1+r_0){\de}_0}. 
 $$
At the beginning of Section 2 we assumed $(2+r_0){\de}_0>1$, hence it is easy 
to see that the proof of (2.6) will be complete if we show
$a_h-a=O(h^{(2+r_0){\de}_0})$. However writing (7.3) with $\a ''=0$ and 
using $\int y^{\a}{\g}(y)\, dy\, =0$ when $1\le |\a |\le 2$ we can replace 
$a(y)-a(x)$ by 
  $$ 
 a(y)-\sum_{|\a |\le 2} (y-x)^{\a}  {\p}^{\a}a(x)/\a !=
  $$
   $$ 
 \sum_{|\a |=2} (y- x)^{\a} \int_0^1 ds\, 2(1-s) \, 
 ({\p}^{\a}a(x+s(y-x))-{\p}^{\a}a(x))/\a ! . 
  $$ 
Since the last expression is $O(|y-x|^{2+r_0})$, we obtain
 $$
 |a_h(x)-a(x)|\le C\int_{{\R}^d} |y-x|^{2+r_0}|{\g}^0_h(x-y)|\, dy\, 
 =C'_{r_0}h^{(2+r_0){\de}_0}. 
 $$
The proof of the assertion b) is described in Appendix of [24]. 
\par \bigskip 
{\bf Proof of Lemma 2.2.} Due to (2.7) and the min-max principle, 
 it suffices to show that (1.5) holds with ${\cal N}(P_h^{\pm},E)$ 
 instead of ${\cal N}(A_h,E)$. 
We drop $\pm$ and we observe that it suffices to prove 
 $$
 {\cal N}(P_h,E)=(2\pi h)^{-d}c^h_E +O(h^{-d})
{\cal R}^{\e ,a_0}_E(h). \eqno (7.5) 
 $$ 
where $c^h_E={\rm vol}\, \{ v\in {\R}^{2d}:\, {\hat p}_h(v)\le E\}$. 
 Indeed, due to (2.6) one has 
 $$
 |c^h_E-c_E|\le {\rm vol}\, \{ v\in {\R}^{2d}:\, |a_0(v)-E|\le Ch\} 
 \le C_{\e}{\cal R}^{\e ,a_0}_E(h). 
 $$
Let $g\in C_0^{\infty}(]E-c ;\; E+c[)$, $l_h=g^2\! \circ \! {\hat p}_h$ 
and $L_h=l_h(x,hD)$. Then reasoning as in Section 3 of [24] we find 
$||g^2(P_h)-L_h||_{\rm tr} =O(h^{1-d})$ and combining this estimate with (2.11) 
we have 
  $$
 {\rm tr\;}(g^2{\tilde f}^Z_h)(P_h) \,=
 \int_{{\R}^{2d}} {dv \over (2\pi h)^d}\, (g^2{\tilde f}^Z_h)
 ({\hat p}_h(v ))   
  \, +O(h^{-d}) \sum_{1\le j\le 2} {\cal R}^{\e /2 ,a_0}_{E_j}(h) \eqno (7.6)
 $$
if $Z=[E_1;\; E_2]\subset [E-c;\; E+c]$. Next we observe that
 $$
|({\tilde f}_h^Z-{\1}_Z)(\la )|\le C_N\sum_{1\le j \le 2} 
{\left( { 1+ { |\la -E_j|\over h } } \right) }^{\! \! -N}  \eqno (7.7)
 $$ 
holds for every $N\in \N$, hence reasoning as in Section 4 we obtain
 $$
 (2\pi h)^{-d}\int_{{\R}^{2d}}(g^2({\tilde f}_h^Z-{\1}_Z))\! \circ \! {\hat p}_h 
  =O(h^{-d})\sum_{1\le j\le 2}{\cal R}^{\e /2,a_0}_{E_j}(h).  \eqno (7.8) 
 $$ 
Taking $Z=[E';\; E'+h]\subset [E-c;\; E+c]$ we obtain 
 $$
 (2\pi h)^{-d}\int_{{\R}^{2d}}(g^2{\tilde f}_h^{\, [E';\; E'+h]})\! \circ \! 
 {\hat p}_h =O(h^{-d}){\cal R}^{\e /2,a_0}_{E'}(h). \eqno (7.9) 
 $$ 
Indeed, due to (7.8) it suffices to observe that (7.9) holds if 
${\tilde f}_h^{\, [E';\; E'+h]}$ is replaced by
${\1}_{[E';\; E'+h]}$. 
As a consequence of (7.6) and (7.9) we find 
 $$
 {\rm tr\;}(g^2{\tilde f}^{\, [E';\; E'+h]}_h)(P_h)=O(h^{-d})
 {\cal R}^{\e /2,a_0}_{E'}(h). \eqno (7.10)  
 $$
We assume moreover $g\ge 0$ and  
 $g=1$ in a neighbourhood of $E$. Let 
${\tilde g}\in C_0^{\infty}(]-\infty ;\; E[)$ 
satisfy ${\tilde g}+g^2=1$ on 
$[\min \{ \inf {\hat p}_h,\, \inf \si (P_h)\} ;\; E]$. Then
 $$
 c^h_E=\int_{{\R}^{2d}}{\tilde g}\! \circ \! {\hat p}_h\, +\int_{{\R}^{2d}}
 (g^2{\1}_{[E-c;\; E]})\! \circ \! {\hat p}_h ,\eqno (7.11) 
 $$
 $$
 {\cal N}(P,E)={\rm tr\, }{\tilde g}(P_h)\, +\, 
 {\rm tr\,}(g^2{\1}_{[E-c;\; E]})(P_h). \eqno (7.12) 
 $$
However reasoning as in Section 3 of [24] we obtain 
 $$
 {\rm tr\, }{\tilde g}(P_h)\, = \, (2\pi h)^{-d}\int_{{\R}^{2d}}  
{\tilde g}\! \circ \! {\hat p}_h\T +\, O(h^{1-d}),   \eqno (7.13)
 $$
hence in order to obtain (7.5) it suffices to show
 $$
{\rm tr\;}(g^2{\1}_{[E-c;\; E]})(P_h)=(2\pi h)^{-d}\int_{{\R}^{2d}}
 (g^2{\1}_{[E-c;\; E]})\! \circ \! {\hat p}_h\, +O(h^{-d})
 {\cal R}^{\e ,a_0}_E(h)  \eqno (7.14)  
 $$ 
and due to (7.6), (7.8), it is clear that (7.14) follows from   
 $$
 {\rm tr\;}(g^2({\tilde f}_h^{\, [E-c;\; E]}-{\1}_{[E-c;\; E]}))(P_h)=O(h^{-d})
{\cal R}^{\e ,a_0}_{E}(h). \eqno (7.15) 
 $$
To begin the proof of (7.15) we observe that modulo $O(h^{\infty})$ 
we can  estimate the left hand side of (7.15) by  
  $$
 {\rm tr\;}g^2(P_h){\left( { 1+ { |P_h-E|\over h } } \right) }^{\! \! -N}
 \le \sum_{k\in \Z}{2\, {\rm tr\,}(g^2{\1}_{[E+kh;\; E+(k+1)h]})(P_h)
 \over (1+ \min \{ |k|,\, |k+1|\})^N} \eqno (7.16) 
 $$
due to (7.7) and clearly the contribution of $\sum_{|k|\ge h^{-\e /2}}$ is 
$O(h^{\e N/2-d})$. 
\par 
Next we observe that ${\tilde \g}_1>0$ allows us to find a constant $C_0>0$ such that 
${\1}_{[E';\; E'+h]}\le C_0{\tilde f}^{\, [E';\; E'+h]}_h$, hence the contribution 
of $\sum_{|k|< h^{-\e /2}}$ in the right hand side of (7.16) can be estimated by 
 $$
 C\sup_{|k|< h^{-\e /2}}{\rm tr\, }(g^2{\tilde f}^{\, [E+kh;\; E+(k+1)h]}_h)(P_h)
\le C_{\e}h^{-d}\sup_{{|k|<h^{-\e /2}}}{\cal R}^{\e /2,a_0}_{E+kh}(h) \eqno (7.17)
 $$
due to (7.10). It is clear that (7.17) can be estimated by 
$O(h^{-d}){\cal R}^{\e ,a_0}_E(h)$. $\triangle$


\end{document}